\documentclass[11pt]{article}

\usepackage{amsmath}
\usepackage{amssymb}
\usepackage{amsfonts}
\usepackage{amsthm}

\usepackage[utf8]{inputenc}

\usepackage{graphicx}
\usepackage[autostyle]{csquotes}
\usepackage{fullpage}
\usepackage{slashed}
\usepackage{mathtools}
\usepackage{indentfirst}
\usepackage{bm}
\usepackage{mathrsfs} 
\usepackage{authblk}
\usepackage{tikz}
\usepackage{float}

\allowdisplaybreaks

\newtheorem{theorem}{Theorem}[section]

\newtheorem{conjecture}[theorem]{Conjecture}

\newcommand{\intg}{\mathbb{Z}}
\newcommand{\rational}{\mathbb{Q}}

\newcommand{\complex}{\mathbb{C}}

\newcommand{\rarw}{\rightarrow}

\newcommand{\lb}{\left(}
\newcommand{\rb}{\right)}
\newcommand{\lsb}{\left[}
\newcommand{\rsb}{\right]}
\newcommand{\lac}{\left\{}
\newcommand{\rac}{\right\}}

\newcommand{\ptl}{\partial}

\newcommand{\RN}[1]{%
  \textup{\uppercase\expandafter{\romannumeral#1}}%
}

\title{A Connect Sum Formula for the BPS Invariant of Knot Complements}
\author{John Chae}
\affil{Department of Physics and QMAP, UC Davis, 1 Shields Ave, Davis, CA, 95616, USA \\ {yjchae@ucdavis.edu}}
\date{}  

\begin{document}

\maketitle

\begin{abstract}
A connect sum formula for the two variable series invariant of a complement of knot is proposed. We provide two kinds of numerical evidence for the proposed formula by examining various torus knots. 

\end{abstract}

\tableofcontents

\section{Introduction}

Inspired by a prediction for a categorification of the Witten-Reshitikhin-Turaev invariant of a closed oriented 3-manifold~\cite{W, RT1} in \cite{GPPV, GPV}, a two variable series invariant $F_K(x,q)$ for a complement $M^{3}_K$ of a knot $K$  was introduced in \cite{GM}. Although its rigorous definition is yet to be found, it possesses various properties such as the Dehn surgery formula and the gluing formula. This knot invariant $F_K$ takes the form\footnote{Implicitly, there is a choice of group; originally, the group used is ${\rm SU}(2)$.}
\begin{gather}
F_K(x,q)= \frac{1}{2} \sum_{\substack{m \geq 1 \\ m \ \text{odd}}}^{\infty} \big(x^{m/2}-x^{-m/2}\big)f_{m}(q) \in \frac{1}{2^{c}} q^{\Delta} \mathbb{Z}\big[x^{\pm 1/2}\big]\big[\big[q^{\pm 1}\big]\big],
\end{gather}
where $f_{m}(q)$ are Laurent series with integer coefficients\footnote{They can be polynomials for monic Alexander polynomial of $K$ (See Section 3.2)}, $c \in \mathbb{Z}_{+}$ and $\Delta \in \mathbb{Q}$. Moreover, $x$-variable is associated to the relative ${\rm Spin}^c \big(M^{3}_K,\ptl M^{3}_K \big)$-structures, which is affinely isomorphic to $H^2\big(M^{3}_K, \ptl M^{3}_K ; \mathbb{Z}\big) \cong H_1\big(M^{3}_K;\mathbb{Z}\big)$; it has an infinite order, which is reflected as a series in $F_K$. The rational constant $\Delta$ was investigated in \cite{GPP}, which elucidated its intimate connection to the d-invariant (or the correction term) in certain versions of the Heegaard Floer homology ($HF^{\pm}$) for rational homology spheres. The physical interpretation of the integer coefficients in $f_{m}(q)$ are number of BPS states of 3d $\mathcal{N}=2$ supersymmetric quantum field theory on $M^{3}_{K}$ together with boundary conditions on $\ptl M^{3}_{K}$. One of various properties of $F_K$ is that it  was conjectured to possess similar characteristics as the colored Jones polynomial, for example, q-holomorphicity~\cite{GL}.
\noindent\begin{conjecture}[{\cite[Conjecture 1.6]{GM}}] For any knot $K \subset$ $S^3$, the normalized series $f_{K}(x,q)$ satisfies a linear recursion relation generated by the quantum A-polynomial of $K$ $\hat{A}_K(q,\hat{x},\hat{y})$:
\begin{gather}
\hat{A}_{K}(q, \hat{x},\hat{y}) f_{K}(x,q) = 0,
\end{gather}
where $f_{K}:=F_{K}(x,q)/\big(x^{1/2}-x^{-1/2}\big)$. 
\end{conjecture}
\noindent The actions of $\hat{x}$ and $\hat{y}$ are
$$
\hat{x} f_{K}(x,q)= x f_{K}(x,q) \qquad \hat{y}f_{K}(x,q)= f_{K}(xq,q).
$$
\newline
This property was used to compute $F_K$ for the figure eight knot $4_1$ in \cite{GM} and was verified for $m(5_2)$ in \cite{P2}. Moreover, the same method was applied to find $F_K$ for a cabling of $4_1$~\cite{C2}. 
\newline

\noindent\textbf{Acknowledgments.} I would like to thank Carsten Schneider and Pavel Putrov for helpful explanations. I am grateful to Sergei Gukov for valuable suggestions on the draft of this paper.

%
%
%

\section{A connected sum formula}

We propose a connect sum formula for $F_K$. 

\noindent\begin{conjecture} For any two knots $K_1$ and $K_2$ in $Y=\intg HS^3$, $F_K(x,q)$ of their connect sum $K_1 \, \# \, K_2$ is
\begin{equation}
F_{K_1 \, \# \, K_2} (x,q) = \frac{F_{K_1}(x,q) F_{K_2}(x,q)}{x^{1/2} - x^{-1/2}}\, \in \frac{1}{2^{c}} q^{\Delta} \mathbb{Z}\big[x^{\pm 1/2}\big]\big[\big[q^{\pm 1}\big]\big],
\end{equation}
where $c \in \intg_{+}$ and $\Delta \in \rational$.
\end{conjecture}

\section{Quantum torus and recursion ideal}

Let $\mathcal{T}$ be a quantum torus
$$
\mathcal{T} := \complex[t^{\pm 1}]\left\langle M^{\pm 1}, L^{\pm 1}\right\rangle / (LM- t^2 ML).
$$
The generators of the noncommutative ring $\mathcal{T}$ acts on a set of discrete functions, which are colored Jones polynomials $J_{K,n} \in \intg[t^{\pm 1}]$ in our context, as
$$
M J_{K,n}= t^{2n} J_{K,n} \qquad L J_{K,n}=  J_{K,n+1}.
$$
The recursion(annihilator) ideal $\mathcal{A}_{K}$ of $J_{K,n}$ is the left ideal $\mathcal{A}_{K}$ in $\mathcal{T}$ consisting of operators that annihilates $J_{K,n}$:
$$
\mathcal{A}_{J_{K,n}} : = \lac \alpha_{K} \in \mathcal{T}\, |\, \alpha_{K} J_{K,n} = 0 \rac.
$$ 
As it turns out that $\mathcal{A}_{K}$ is not a principal ideal in general. However, by adding inverse polynomials of t and M to $\mathcal{T}$~\cite{G},
we obtain a principal ideal domain $\tilde{\mathcal{T}}$
$$
\tilde{\mathcal{T}} : = \lac \sum_{j \in \intg} a_{j}(M) L^{j} \Big|\, a_{j}(M) \in \complex[t^{\pm 1}](M),\, a_{j}= \text{almost always}\quad 0 \rac
$$
Using $\tilde{\mathcal{T}}$ we get a principal ideal $\tilde{\mathcal{A}_{K}}:= \tilde{\mathcal{T}}\mathcal{A}_{K}$ generated by a single polynomial $\hat{A}_{K}$
$$
\hat{A}_{K}(t,M,L)= \sum_{j=0}^{d} a_{j}(t,M)L^{j}.
$$
This $\hat{A}_{K}$ polynomial is a noncommutative deformation of a classical A-polynomial of a knot~\cite{CCGLS} (see also \cite{CL}). Alternative approaches to obtain $\hat{A}_K(t,M,L)$ are by quantizing the classical A-polynomial curve using the twisted Alexander polynomial or applying the topological recursion technique~\cite{GS}. The AJ conjecture states that the classical polynomial can be obtained from its quantum version by setting $t=-1$ (up to an overall rational function of $M$)~\cite{G, Gukov}.

\section{Recursion relations}

We provide evidence for (3) using the q-holonomic property (2) of $F_K$ for connected sums of torus knots. For right handed torus knots $T(s,t)\quad 2 \leq s < t\quad gcd(s,t)=1$, their $F_K$ were computed~\cite{GM}:
\begin{equation}
F_K(x,q)= \frac{1}{2} \sum_{\substack{m \geq 1 \\ m \ \text{odd}}}^{\infty} \ \epsilon_m \big(x^{m/2}-x^{-m/2}\big) q^{\frac{m^2-(st-s-t)^2}{4st}}
\end{equation}
$$
\epsilon_m = \begin{cases}
-1,\quad  m \equiv st+s+t\quad \text{or}\quad st-s-t\quad \text{mod}\, 2st\\
+1,\quad  m \equiv st+s-t\quad \text{or}\quad st-s+t\quad \text{mod}\, 2st\\
0,\quad \text{otherwise.}
\end{cases}
$$
For the left handed torus knots $T(s,-t)$, their coefficient functions can be obtained from (4) by $f_m (q^{-1})$ and $F_{T(s,-t)} \in 2^{-c} q^{\Delta} \mathbb{Z}\big[q^{\pm 1}\big]\big[\big[x^{\pm 1/2}\big]\big]$. 
\newline

\noindent In the following examples, we used \cite{KK} to obtain the quantum A-polynomials for the connected sum of knots.
\newline

\noindent \underline{$K=T(2,3)\, \# \, T(2,3)$}\quad The minimal degree homogeneous recursion relation for $K$ is
\begin{equation}
r_0\, F_{T(2,3)}(x,q)^2 + r_1\, F_{T(2,3)}(xq,q)^2 + r_2\, F_{T(2,3)}(xq^2,q)^2 + r_3\, F_{T(2,3)}(xq^3,q)^2 = 0
\end{equation}
\begin{align*}
r_0 & = -q+q \left(q^3+q^5\right) x^2+q^6 x^3-q^9 x^4-2 q^{11} x^5+q^{16} x^7\\
r_1 & = 1+\left(-2 q+q^3-q^5\right) x^2-2 q^2 x^3+\left(q^2-2 q^4+3 q^6-q^8\right) x^4+\left(2 q^3+2 q^7\right)x^5\\
& +\left(q^4 +q^5-2 q^7+3 q^9-q^{11} \right) x^6-2 q^8 x^7+\left(-q^7-q^9-2q^{10}+q^{12}-q^{14}\right) x^8 -q^9 x^9\\
& +\left(q^{12}  +q^{15}  \phantom{1} \right) x^{10}+2 q^{14} x^{11}-q^{19}x^{13}\\
r_2 & = q^4 x^3-2 q^5 x^5+q^4 \left(-q^2-q^5\right) x^6+q^6 x^7+q^4 \left(q^3+q^5+2 q^6-q^8+q^{10}\right) x^8+2q^{11} x^9\\
& +q^4 \left(-q^6-q^7+2 q^9-3 q^{11}+q^{13}\right) x^{10}+q^4 \left(-2 q^8-2 q^{12}\right) x^{11}+q^4 \left(-q^{10}+2 q^{12}-3 q^{14}+q^{16}  \right) x^{12}\\
& +2 q^{17} x^{13}+q^4 \left(2q^{15}-q^{17}+q^{19}\right) x^{14}-q^{24} x^{16}\\
r_3 & = -q^{19} x^9+2 q^{20} x^{11}+q^{21} x^{12}-q^{21} x^{13}+q^{19} \left(-q^3-q^5\right) x^{14}+q^{25}x^{16}\\
\end{align*}
For instance, at x-order, the above four terms in the same order are
\begin{align*}
T_0 & =-q^3+3 q^6+q^8-2 q^{11}+2 q^{12}-2 q^{13}+8 q^{15}-2 q^{17}-2 q^{18}-4 q^{20}-2q^{21}+O(q^{23})\\
T_1 & =-4 q^2-4 q^3+6 q^4-14 q^6-4 q^7+4 q^8-6 q^9+q^{10}+3 q^{11}-10 q^{12}+q^{13}+12q^{14}- O( q^{15})\\
T_2 & = \frac{3}{q} -1+9 q+8 q^2+3 q^3-4 q^4+13 q^6+8 q^7-5 q^8+12 q^9-q^{10}+q^{11}+10q^{12}+q^{13}-O(q^{14})\\
T_3 & = -\frac{3}{q} + 1 -9 q-4 q^2+2 q^3-2 q^4-2 q^6-4 q^7-6 q^9-2 q^{11}-2 q^{12}-2 q^{17}+4q^{19}-O(q^{21})\\
\end{align*}
The figure below shows that as the upper bound of the summation in (1) increases, the minimum power of q-term that survived increases. This indicates that the desired cancellations occur. 
\begin{figure}[h!]
\begin{center}
\includegraphics[scale=1]{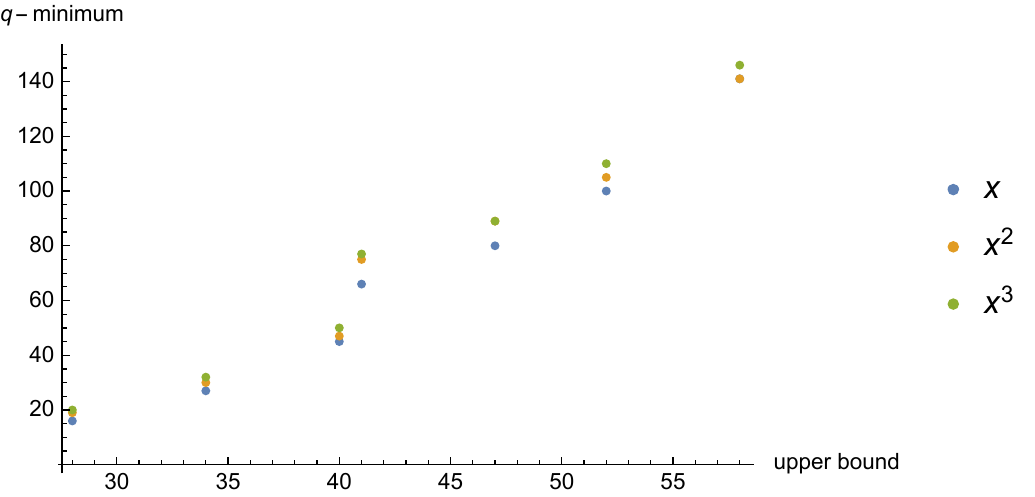}
\caption{Minimum powers of q-terms that survived in (5) for the powers of x shown. The upper bound corresponds to the maximum value among the upper bounds in the summations in $F_{T(2,3)  \, \# \, T(2,3)} (xq^j,q),\, j=0,\cdots ,3$\, (see Appendix A for the plots of other powers of x).}
\end{center}
\end{figure}
\newline

\noindent \underline{$K=T(2,3)\, \# \, T(2,5)$}\quad The minimal degree homogeneous recursion relation for $K$ is
$$
t_0\, F_{T(2,3)}(x,q)\, F_{T(2,5)}(x,q) + t_1\, F_{T(2,3)}(xq,q)\, F_{T(2,5)}(xq,q) + t_2\, F_{T(2,3)}(xq^2,q)\,F_{T(2,5)}(xq^2,q)
$$
\begin{equation}
\hspace{-3cm} + t_3\, F_{T(2,3)}(xq^3,q)\, F_{T(2,5)}(xq^3,q) + t_4\, F_{T(2,3)}(xq^4,q)\, F_{T(2,5)}(xq^4,q) = 0
\end{equation}
The coefficient functions $t_i (x,q) \in \intg[x,q]$ are recorded in \cite{C1}. At x-order, the above five terms in the same order are
\begin{align*}
R_0 & = -q^5+q^6-q^7+q^8-2 q^9+5 q^{10}-q^{11}+q^{12}-q^{13}+4 q^{14}+2 q^{15}+4 q^{16}-3 q^{17}- O( q^{19})\\
R_1 & = -4 q^4-q^5+7 q^6+3 q^7-7 q^8-11 q^9-8 q^{10}+10 q^{11}+14 q^{12}-5 q^{13}-2 q^{14}-7 q^{15}- O(q^{16})\\
R_2 & = -\frac{2}{q^2}+\frac{4}{q} -10 +12 q-5 q^2-3 q^3+2 q^4+3 q^5-q^6+7 q^7-11 q^8-14 q^9+4 q^{10}-O( q^{11})\\
R_3 & = -\frac{1}{q^7}+\frac{1}{q^6}-\frac{4}{q^5}+\frac{6}{q^4}+\frac{1}{q^3}-\frac{5}{q^2}+\frac{3}{q}+11 -13q+2 q^2+6 q^3+3 q^4-6 q^5-9 q^6-O( q^7)\\
R_4 & = \frac{1}{q^7}-\frac{1}{q^6}+\frac{4}{q^5}-\frac{6}{q^4}-\frac{1}{q^3}+\frac{7}{q^2}-\frac{7}{q}-1+q+3q^2-3 q^3-q^4+5 q^5+2 q^6-q^7+q^9-O( q^{10})\\
\end{align*}

\begin{figure}[h!]
\begin{center}
\includegraphics[scale=1]{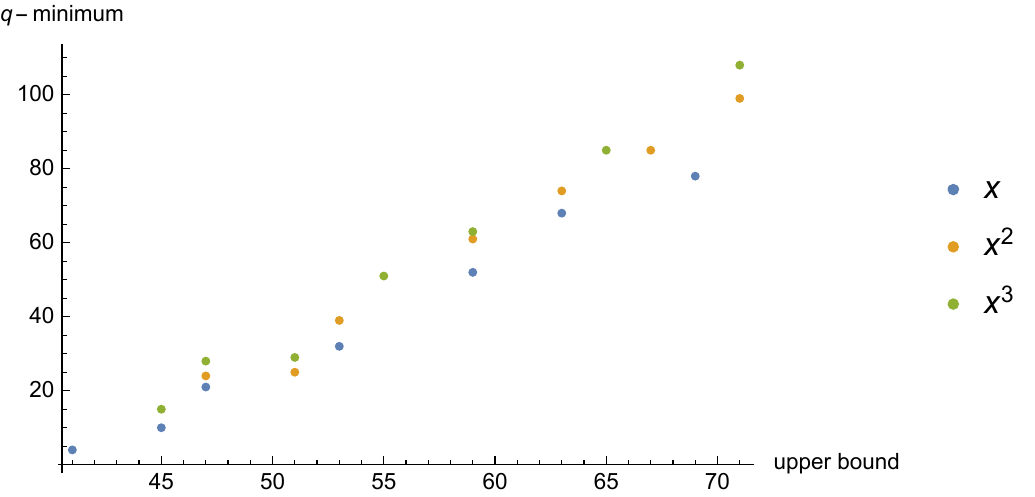}
\caption{Minimum powers of q-terms that survived in (6) for the powers of x shown. The upper bound corresponds to the maximum value among the upper bounds in the summations in $F_{T(2,3)  \, \# \, T(2,5)} (xq^j,q),\, j=0,\cdots ,4$\, (see Appendix A for the plots of other powers of x).}
\end{center}
\end{figure}

\noindent \underline{$K=T(2,3)\, \# \, T(3,5)$}\quad The minimal degree homogeneous recursion relation for $K$ is
$$
h_0\, F_{T(2,3)}(x,q)\, F_{T(2,5)}(x,q) + h_1\, F_{T(2,3)}(xq,q)\, F_{T(2,5)}(xq,q) + h_2\, F_{T(2,3)}(xq^2,q)\, F_{T(2,5)}(xq^2,q)
$$
\begin{equation*}
\hspace{-4cm} + h_3\, F_{T(2,3)}(xq^3,q)\, F_{T(2,5)}(xq^3,q) + h_4\, F_{T(2,3)}(xq^4,q)\, F_{T(2,5)}(xq^4,q) 
\end{equation*}
\begin{equation}
\hspace{-3.5cm}  + h_5\, F_{T(2,3)}(xq^5,q)\, F_{T(2,5)}(xq^5,q) + h_6\, F_{T(2,3)}(xq^6,q)\, F_{T(2,5)}(xq^6,q)  = 0
\end{equation}
The coefficient functions $h_i (x,q)\in \intg[x,q]$ are listed in \cite{C1}.  At $x^{0}$ order, the above seven terms in the same order are
\begin{align*}
W_0 & = 5 q^{304}+9 q^{305}+15 q^{306}+32 q^{307}+48 q^{308}+64 q^{309}+79 q^{310}+91 q^{311}+ O(q^{312})\\
W_1 & =q^{295}+2 q^{296}-5 q^{297}-18 q^{298}-36 q^{299}-66 q^{300}-104 q^{301}-155 q^{302}- O(q^{303})\\
W_2 & = -q^{286}+q^{288}-4 q^{289}-10 q^{290}-11 q^{291}+9 q^{292}+41 q^{293}+78 q^{294}+O(q^{295})\\
W_3 & = 2 q^{289}+9 q^{290}+4 q^{291}-9 q^{292}-16 q^{293}-35 q^{294}-51 q^{295}-54 q^{296}- O( q^{297})\\
W_4 & = q^{275}+2 q^{276}-q^{277}-8 q^{278}+4 q^{279}+13 q^{280}+19 q^{281}+32 q^{282}+46 q^{283}+O(q^{284})\\
W_5 & = -2 q^{289}-9 q^{290}-4 q^{291}+9 q^{292}+16 q^{293}+35 q^{294}+50 q^{295}+52 q^{296}+50 q^{297}+O( q^{298})\\
W_6 & = -q^{275}-2 q^{276}+q^{277}+8 q^{278}-4 q^{279}-13 q^{280}-19 q^{281}-32 q^{282}-46 q^{283}-45q^{284}-O( q^{285})\\
\end{align*}

\begin{figure}[h!]
\begin{center}
\includegraphics[scale=1]{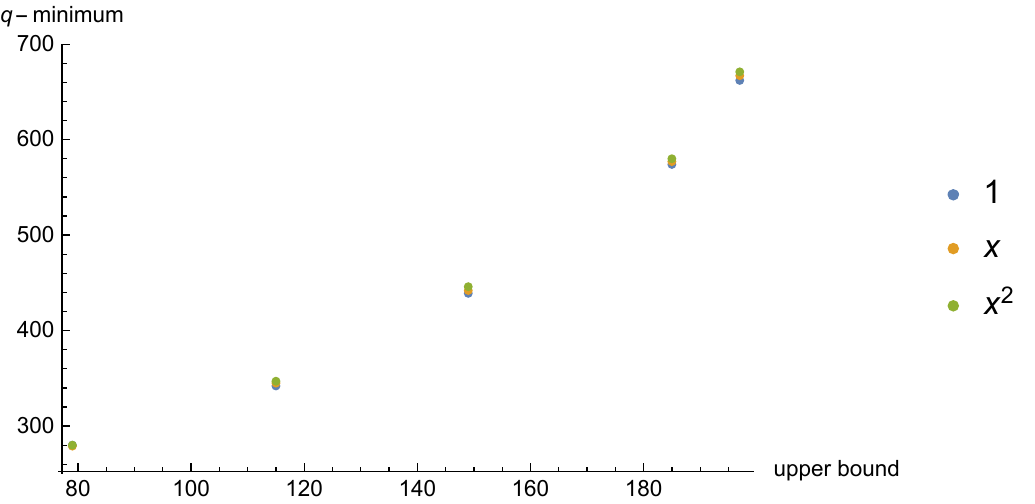}
\caption{Minimum powers of q-terms that survived in (7) for the powers of x shown. The three dots are nearly overlapping. The upper bound corresponds to the maximum value among the upper bounds in the summations in $F_{T(2,3)  \, \# \, T(3,5)} (xq^j,q),\, j=0,\cdots ,6$\, (see Appendix A for the plots of other powers of x).}
\end{center}
\end{figure}

\noindent \underline{$K=T(2,3)\, \# \, T(2,-3)$}\quad The minimal degree homogeneous recursion relation for $K$ is
$$
b_0\, F_{T(2,3)}(x,q)\, F_{T(2,-3)}(x,q) + b_1\, F_{T(2,3)}(xq,q)\, F_{T(2,-3)}(xq,q) + b_2\, F_{T(2,3)}(xq^2,q)\, F_{T(2,-3)}(xq^2,q)
$$
\begin{equation}
\hspace{-4cm} + b_3\, F_{T(2,3)}(xq^3,q)\, F_{T(2,-3)}(xq^3,q)+ b_4\, F_{T(2,3)}(xq^4,q)\, F_{T(2,-3)}(xq^4,q)=0,
\end{equation}
\newline
\noindent where the coefficient functions $b_i (x,q)\in \intg[x,q]$ are listed in \cite{C1}. For this composite knot, the cancellation in (8) is subtle compared to the connected sums of the right handed torus knots since the coefficient function $f_m(q)$ of the left handed torus knots have the form $q^{-r},\, r \in \intg_{+}$. Specifically, arbitrary high and low powers of q from $F_{T(2,3)}$ and $F_{T(2,-3)}$, respectively, which appear for large values of the upper bound of the summations in $F_{T(2,\pm 3)}$, can combine to yield $O(1)$-powers of q that is required for cancellations. Desired cancellations become evident when we group the terms in (8) in powers of q and observe cancellations among x terms. It turns out that for some powers of q such as $q$ (Figure 20)  and $q^{500}$ (Figure 27), cancellations do not occur in $x^p$ or $x^{-p},\quad p\in \intg_{+}$ when the upper bound is not high enough. Furthermore, another gap can be created for some powers of q when the upper bound is high enough.  Therefore, we scrutinized the growth of width of gaps in x-terms as the upper bound is increased for various powers of q.
\newline

For example, when the upper bound of the summation is 325, a subset of x-terms at $q^{100}$ in (8) are
\begin{align*}
(8) & \supset \frac{76}{x^{281}}+\frac{49}{x^{280}}-\frac{118}{x^{279}}-\frac{21}{x^{278}} +\frac{51}{x^{277}}+\frac{26}{x^{276}} -\frac{11}{x^{275}}-\frac{34}{x^{274}}+\frac{14}{x^{273}}+\frac{14}{x^{272}}-\frac{3}{x^{271}}-\frac{3}{x^{270}}\\
& -2 x^2-8 x^3+16 x^4 +91 x^5 -83 x^6-151 x^7+69 x^8+154 x^9-71 x^{10}-15 x^{11}-x^{12}+x^{13}\\
& -7   x^{287}+24 x^{288}-6 x^{289} -33 x^{290}+14 x^{291} +14 x^{292}-12 x^{293}+15 x^{294} +6 x^{295}-24x^{296}\\
&  +x^{297}+12 x^{298}+3 x^{299}-55 x^{300}+28 x^{301} +100 x^{302}-57 x^{303}  -56 x^{304}+48 x^{305}-47 x^{306}+ \cdots
\end{align*}
\noindent There is a gap between $x^{14}$ and $x^{286}$ and there is another gap from $x^0$	to $x^{-269}$. These gaps are due to cancellations as we can see from the five terms in Appendix A. In the figure below, we observe that the gap size widens for $q^{100}$ as the upper bound of the summation is increased.
\begin{figure}[h!]
\begin{center}
\includegraphics[scale=1]{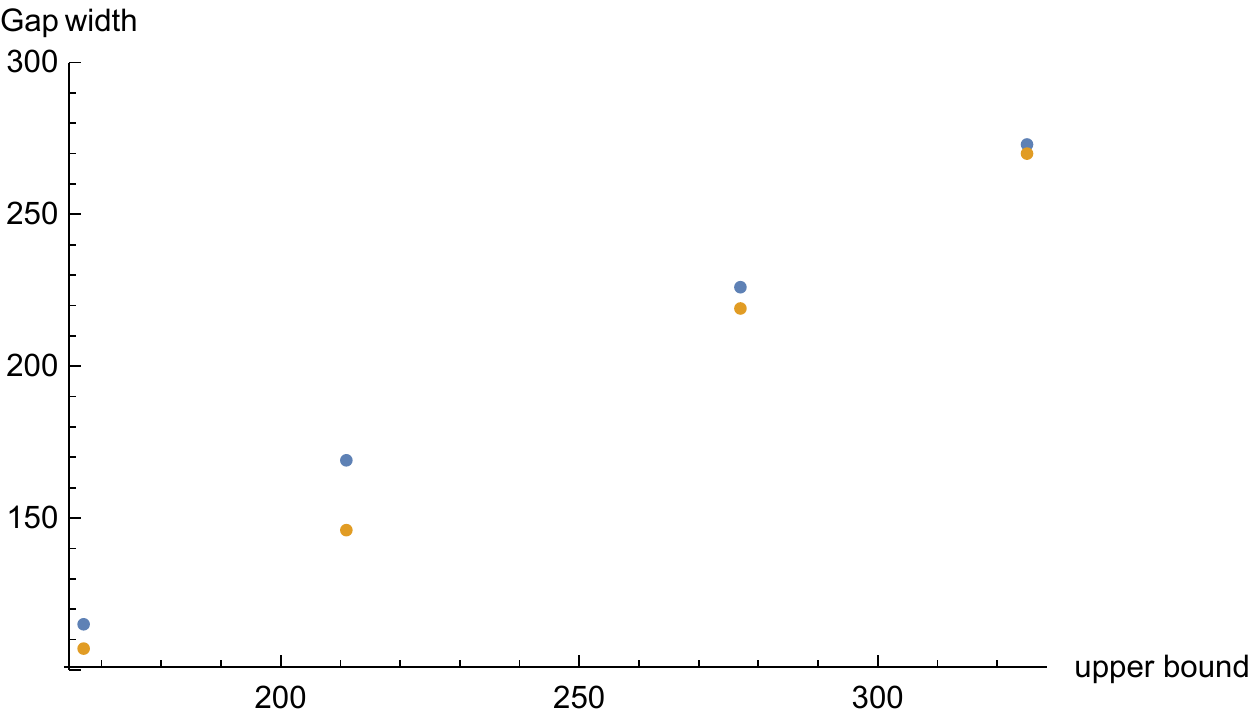}
\caption{At $q^{100}$, the width of the gaps in $x^p$ terms (blue) and in $1/x^p$ terms (orange), $p \in \intg_{+}$ is displayed. The upper bound corresponds to the maximum value among the upper bounds in the summations in $F_{T(2,3)  \, \# \, T(2,-3)} (xq^j,q),\, j=0,\cdots ,4$\, (see Appendix A for the plots of other powers of q).}
\end{center}
\end{figure}
\newline
For lower powers of q such as $q,q^2$ and $q^{-1}$, when the upper bound is 165, cancellations among small positive powers of x occur. This second gap widens as the upper bound is increased. For instance, at $q^{-3}$, $x$ and $x^2$ terms are absent when the upper bound is 165. As it is increased to 187, $x$ to $x^5$ terms are canceled.
\newline
%

\section{Comparison to the analytic results}

In this section, we compare the $SU(2)$ WRT invariant of integral homology spheres at fixed roots of unity obtained analytically and numerically. For the latter method, we utilize the conjectured Dehn surgery formula in \cite{GM}, which relates $F_K$ and  $\hat{Z}$:
\newline

\noindent \textbf{Conjecture 5.1} ({\cite[Conjecture 1.7]{GM}})\, For any $K \subset S^3$ and let $S^3_{p/r}(K)$ be a 3-manifold obtained from Dehn surgery on $K$ along $p/r \in \rational^{\ast}$. Then
	$$
	\hat{Z}_{b}[S^3_{p/r}(K); q] = \pm q^d \, \mathcal{L}^{(b)}_{p/r} \lsb \lb x^{\frac{1}{2r}} - x^{-\frac{1}{2r}} \rb F_{K}(x,q) \rsb \qquad d \in \rational,
	$$
	$$
	\mathcal{L}^{(b)}_{p/r} : x^{u}q^{v} \mapsto \begin{cases}
      q^{-u^2 r/p}q^v & \text{if}\quad ru - b \in p\intg \\
      0 & \text{otherwise}
    \end{cases} 
	$$
	where $\mathcal{L}$ is a $|q|<1$ generalization of the Laplace transform~\cite{BBL}.
\newline

\noindent On analytic side, the integer Dehn surgery formula for the WRT invariant at a primitive $k$-th root of unity is~\cite{BBL,BL}
\begin{equation}
\tau_k [S^3_{p}(K)] = \frac{\sum_{n=1}^{k-1} \, [n]^2 \, q^{p(n^2 -1)/4}\, J_{n}(K)}{\sum_{n=1}^{k-1}\, [n]^2 \, q^{sign(p) (n^2 -1)/4}}\qquad [n]= \frac{q^{n/2}-q^{-n/2}}{q^{1/2}-q^{-1/2}}
\end{equation}
where $J_{n}(K)$ is $sl(2)$ colored Jones polynomial of $K$ and $p \in \intg$ is the surgery slope or equivalently framing of $K$. When $p=-1$, it results in $S^3_{-1}(K)=\intg HS^3$ for any $K$. For this class of manifolds, the decomposition of the $SU(2)$ WRT invariant in terms of $\hat{Z}$ is~\cite{GPPV}
\begin{equation}
Z_{CS} \lsb S^3_{-1}(K); q=e^{\frac{i 2 \pi}{k}} \rsb = \frac{-i}{2\sqrt{2k}}\, \lim_{q \rarw e^{\frac{i 2 \pi}{k}}} \hat{Z}_{0} (q).
\end{equation}
It is simply related to $\tau_k$
$$
Z_{CS} \lsb S^3_{-1}(K); q=e^{\frac{i 2 \pi}{k}} \rsb = \frac{-i(q^{1/2} - q^{-1/2})}{\sqrt{2k}}\, \tau_k [S^3_{-1}(K)].
$$
For the examples below, we display the $sl(2)$ colored Jones polynomial for the torus knot $T(s,t)$,
	$$
	J_{n}(T(s,t);q)= -\frac{q^{-\frac{s t n^2}{4} } q^{\frac{(s-1) (t-1)}{2} }}{q^{\frac{n}{2}}-q^{-\frac{n}{2}}}\, \sum_{r=0}^{stn} \epsilon_{s t n-r}\,  q^{\frac{r^2-(s t-s-t)^2}{4 s t}}\qquad n \in \mathbb{N}
	$$
where $ 2 \leq s < |t|$,\, $gcd(s,t)=1$ and $\epsilon$ is in (4) (The unknot normalization is $J_{n}=1$).
\newline

\noindent \underline{$K=T(2,3)\, \# \, T(2,3)$}: At $k=3$, applying the analytic formula (9) yields 
	$$
	Z_{CS} \lsb S^3_{-1}(K); e^{\frac{i 2 \pi}{3}} \rsb = 0.7071,
	$$
where $J_{n}(K_1 \, \# K_2)= J_{n}(K_1) J_{n} (K_2)$ is used. On the numerical side, after $\hat{Z}$ is obtained from Conjecture 5.1, we truncate the q-power series at a large finite power $N$ of q to find the limiting value of $\hat{Z}_{0} (q)$ as $q$ goes to a root of unity. We choose the truncation power to be $N=20000$ and extract the limiting value of $\hat{Z}_{0} (q)$. The figure below shows that the q-series converges to 
	$$ 
	\lim_{q \rarw e^{\frac{i 2 \pi}{3}}}\, \frac{2}{q^2}\hat{Z}_{0} [S^3_{-1}(K); q] \longrightarrow   -0.0003504774588  - i 6.925958533.
	$$
The overall monomial is introduced for numerical convenience. After substituting the limiting value into (10), we find $Z_{CS} \approx 0.7068717087$, thus it agrees with the above analytical value. 
\begin{figure}[H]
\begin{center}
\includegraphics[scale=0.5]{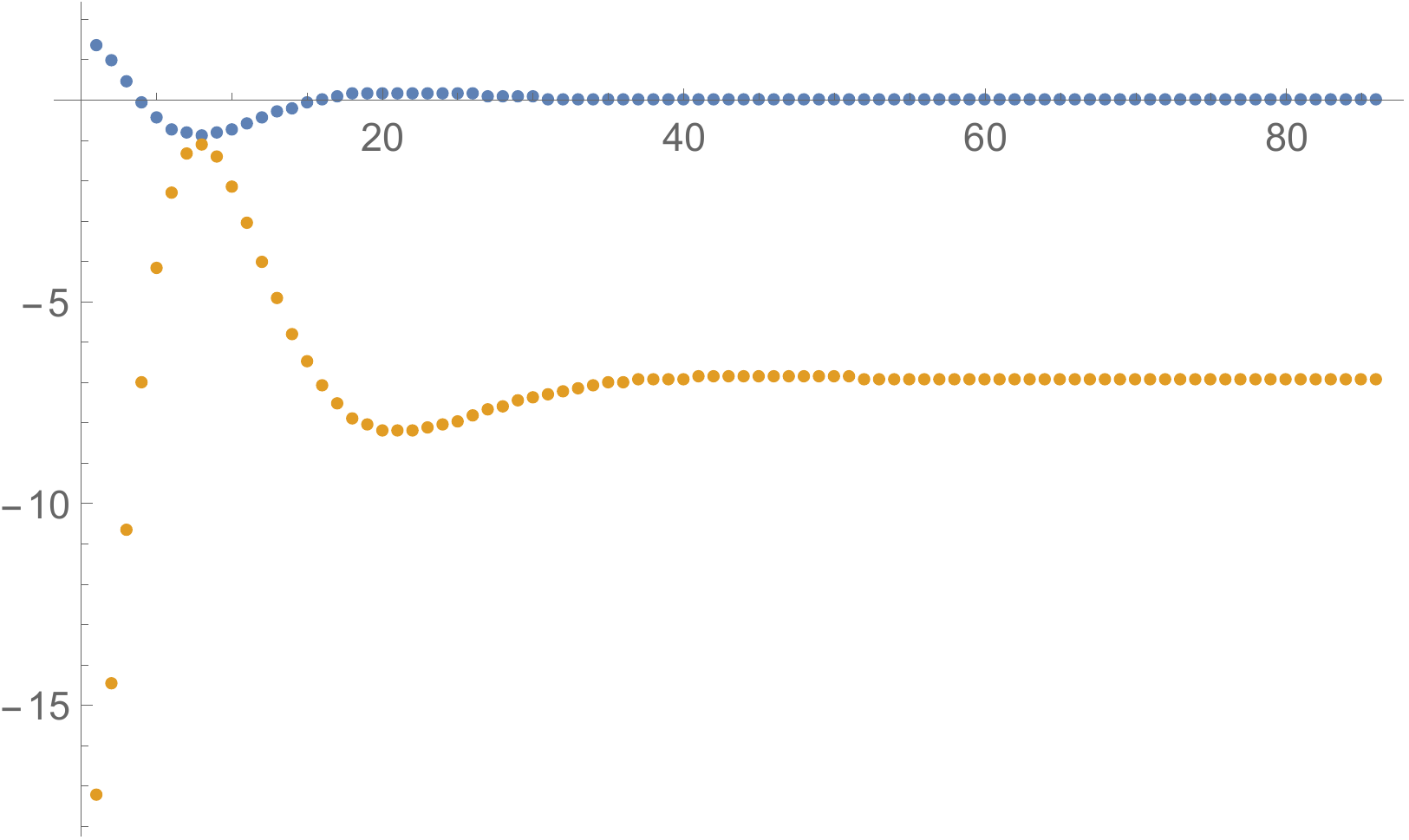}
\caption{The extrapolation of $2\hat{Z}_{0} (q \rarw e^{\frac{i 2 \pi}{3}} )/q^2$ associated with $K$ at $N=20000$. Real part (blue) and imaginary part (orange) of $\hat{Z}_{0}$. }
\end{center}
\end{figure}
\noindent At $k=4$, the analytic formula (9) results in 
	$$
	Z_{CS} \lsb S^3_{-1}(K); e^{\frac{i 2 \pi}{4}} \rsb = 0.5 . 
	$$
As in the previous case, we truncate the q-power series at $N=20000$ and find the limiting value of $\hat{Z}$ as $q$ goes to $i$. 
\begin{figure}[H]
\begin{center}
\includegraphics[scale=0.5]{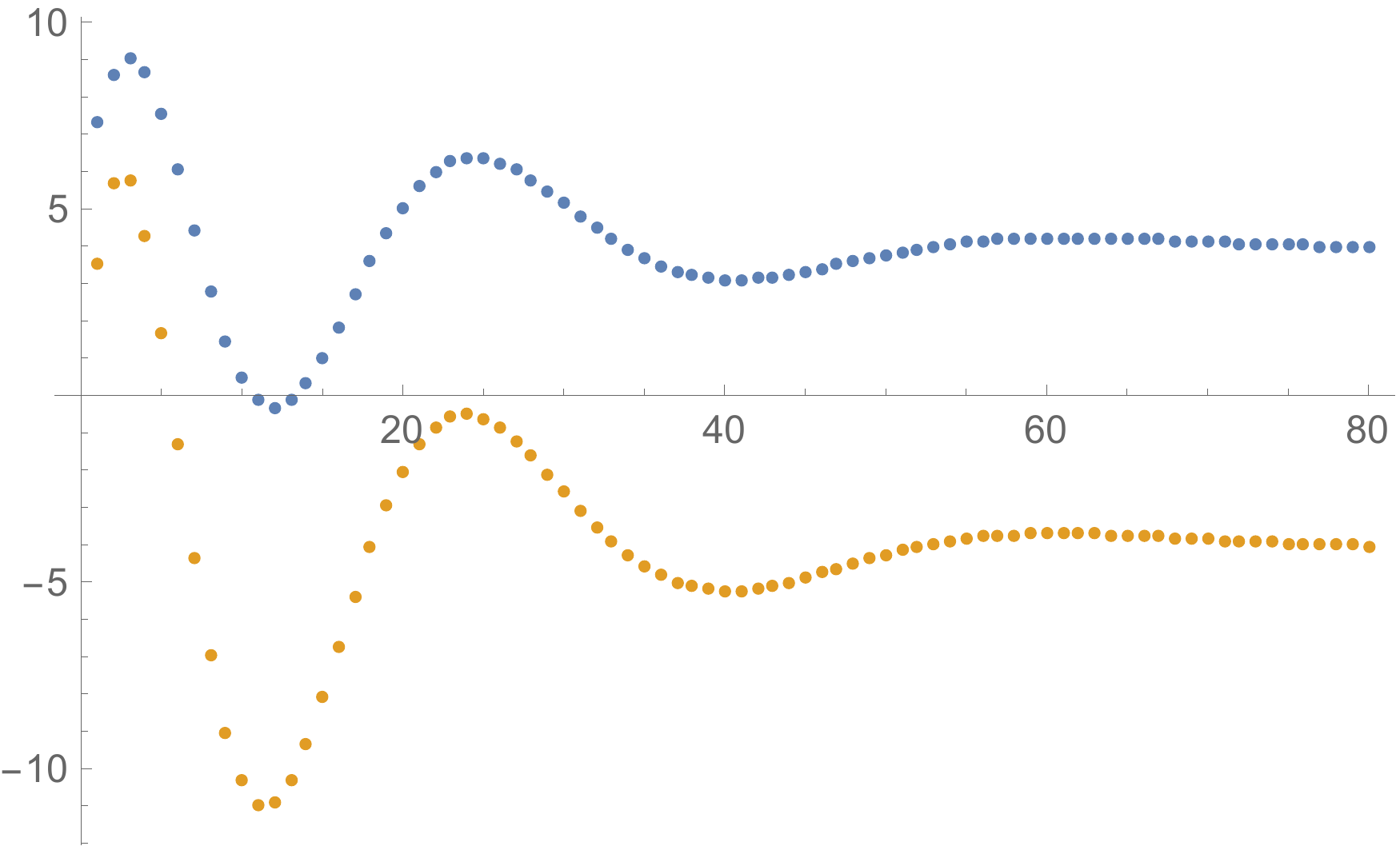}
\caption{The extrapolation of $2\hat{Z}_{0} (q \rarw e^{\frac{i 2 \pi}{4}} )/q^2$ associated with $K$ at $N=20000$. Real part (blue) and imaginary part (orange) of $\hat{Z}_{0}$. }
\end{center}
\end{figure}	
\noindent The q-series approaches to 
	$$ 
	\lim_{q \rarw e^{\frac{i 2 \pi}{4}}}\, \frac{2}{q^2}\hat{Z}_{0} [S^3_{-1}(K); q] \longrightarrow 3.968560094 - i 4.028195455.
	$$
Using (10), $Z_{CS} \approx 0.5 $, which matches with the analytical result. 
\newline
		
\noindent At $k=5$, the analytic formula (9) produces 
	$$
	Z_{CS} \lsb S^3_{-1}(K); e^{\frac{i 2 \pi}{5}} \rsb = -0.3 + i 1.36263. 
	$$
	We truncate the q-power series at $N=30000$ and find the limiting value of $\hat{Z}$. 
\begin{figure}[H]
\begin{center}
\includegraphics[scale=0.5]{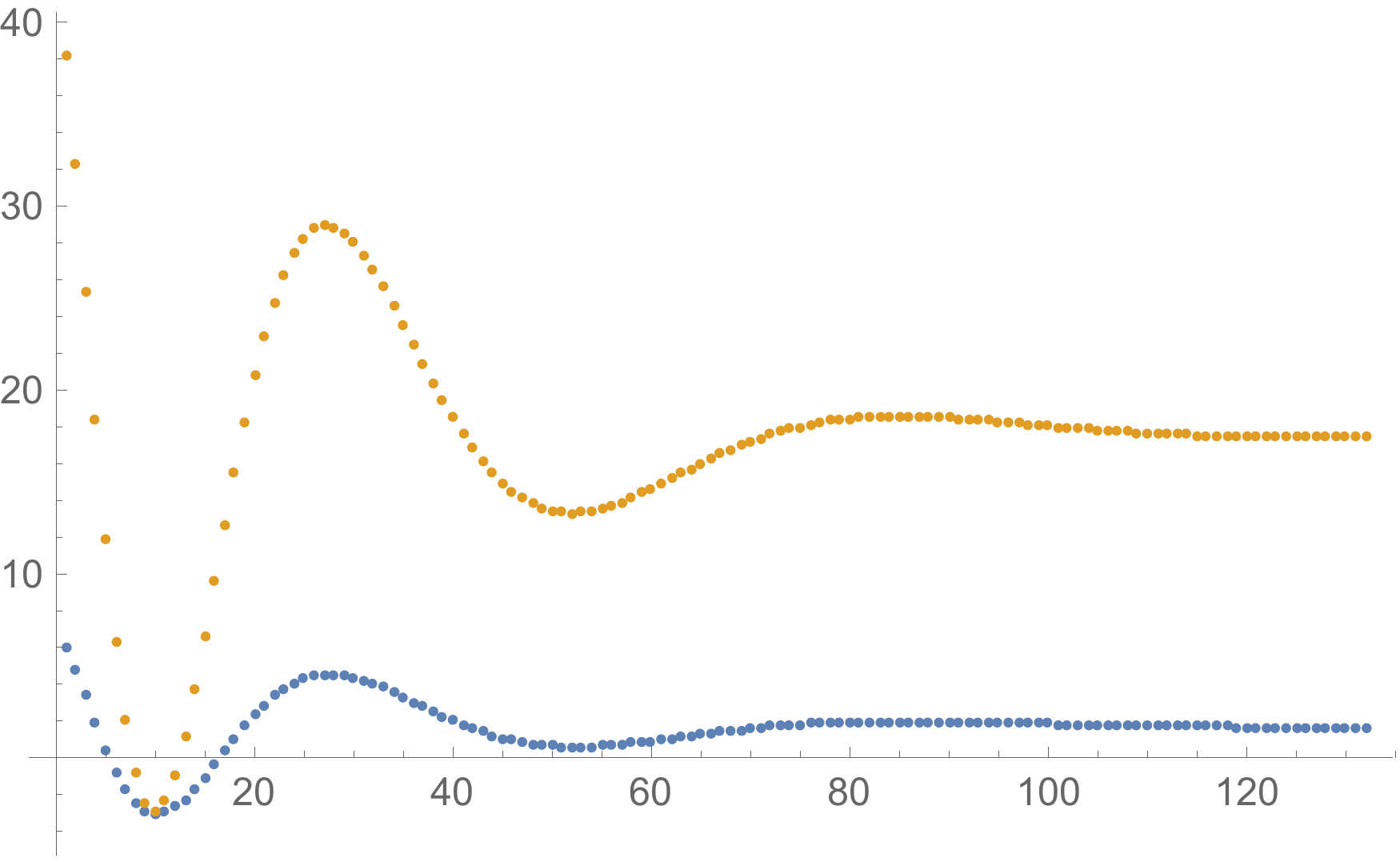}
\caption{The extrapolation of $2\hat{Z}_{0} (q \rarw e^{\frac{i 2 \pi}{5}} )/q^2 $ associated with $K$ at $N=30000$. Real part (blue) and imaginary part (orange) of $\hat{Z}_{0}$. }
\end{center}
\end{figure}	
\noindent The q-series approaches to 
	$$ 
	\lim_{q \rarw e^{\frac{i 2 \pi}{5}}}\, \frac{2}{q^2}\hat{Z}_{0} [S^3_{-1}(K); q] \longrightarrow 1.6675682 + i 17.42149573.
	$$
From (10), $Z_{CS} \approx -0.3 + i 1.35 $, which agrees with the analytical result.
\newline

\noindent \underline{$K=T(2,3)\, \# \, T(2,5)$}: At $k=3$, applying the analytic formula (9) yields 
	$$
	Z_{CS} \lsb S^3_{-1}(K); e^{\frac{i 2 \pi}{3}} \rsb = 0.7071.
	$$
After truncating the q-power series at $N=25000$ and then extracting the limiting value of $\hat{Z}_{0} (q)$ results in Figure 8. It shows that the q-series converges to 
	$$ 
	\lim_{q \rarw e^{\frac{i 2 \pi}{3}}}\, \frac{2}{q^4} \hat{Z}_{0} [S^3_{-1}(K); q] \longrightarrow   5.989718 + i 3.450427632 .
	$$
After substituting it into (10), we find $Z_{CS} \approx 0.705499 - i 0.00068351$.
\begin{figure}[H]
\begin{center}
\includegraphics[scale=0.5]{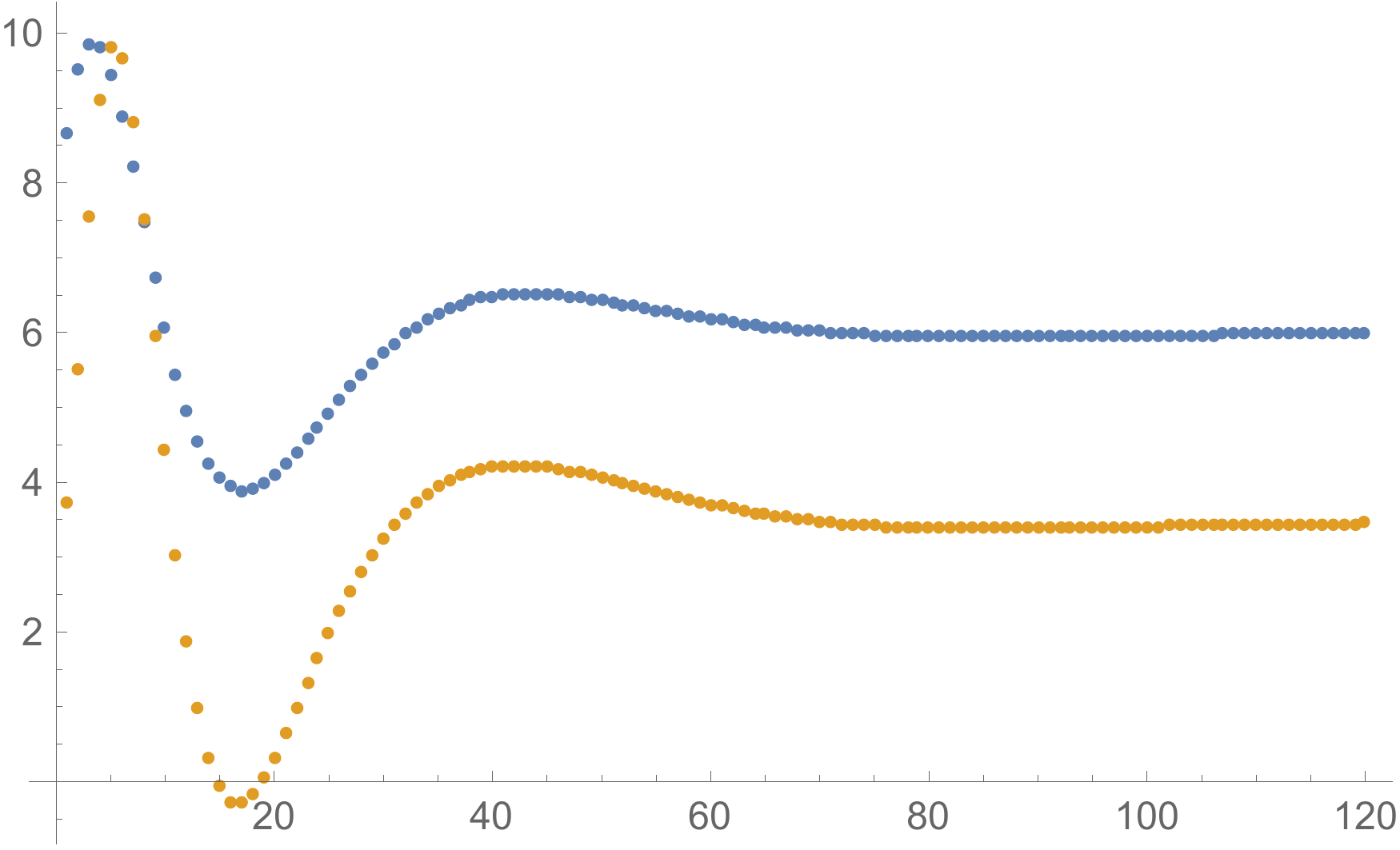}
\caption{The extrapolation of $2\hat{Z}_{0} (q \rarw e^{\frac{i 2 \pi}{3}} )/q^4$ associated with $K$ at $N=25000$. Real part (blue) and imaginary part (orange) of $\hat{Z}_{0}$. }
\end{center}
\end{figure}
\noindent At $k=4$, the analytic formula (9) results in 
	$$
	Z_{CS} \lsb S^3_{-1}(K); e^{\frac{i 2 \pi}{4}} \rsb = 0.5 . 
	$$
As in the above case, we truncate the q-power series at $N=25000$ and find the limiting value of $\hat{Z}$. 
\begin{figure}[H]
\begin{center}
\includegraphics[scale=0.5]{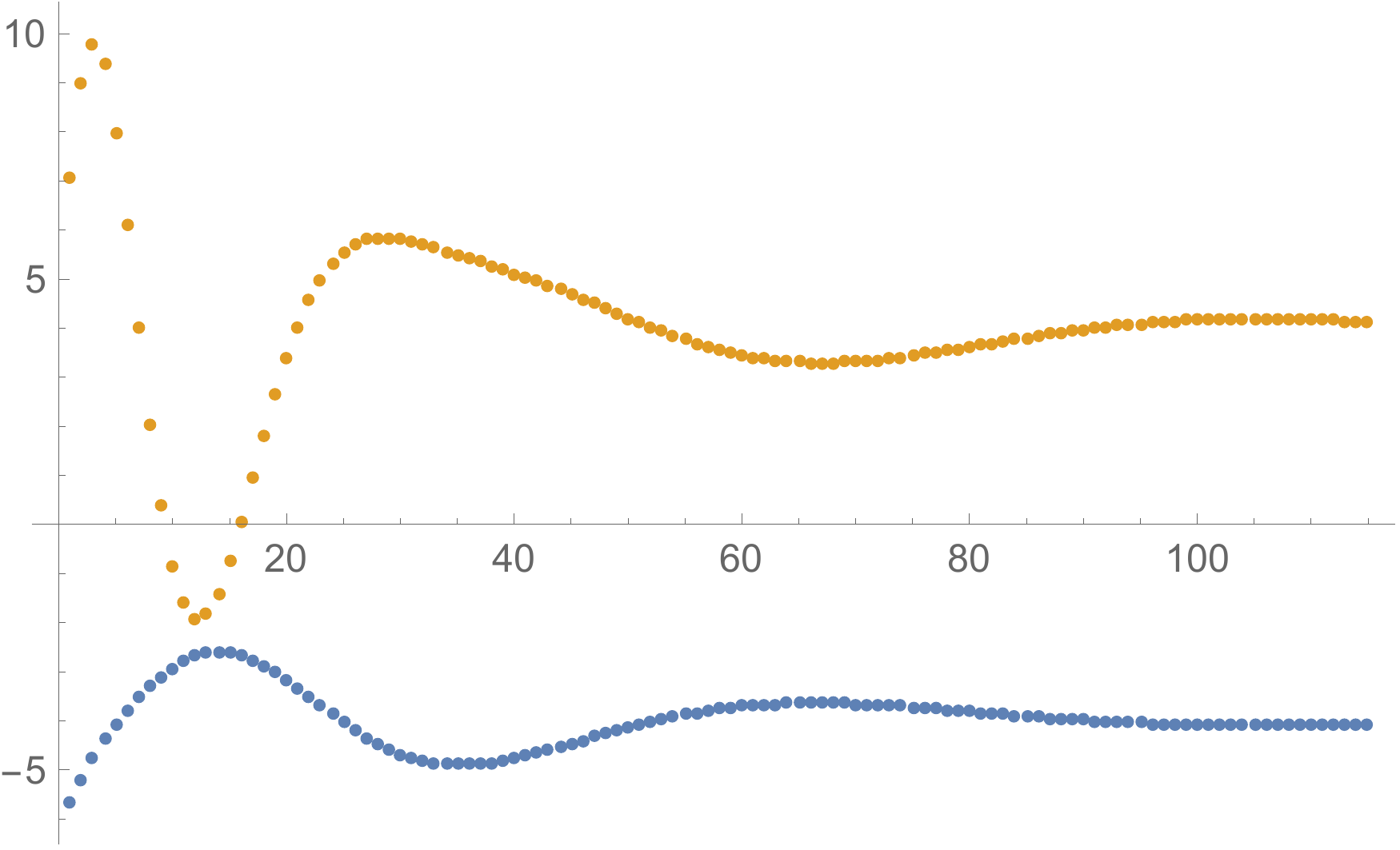}
\caption{The extrapolation of $2\hat{Z}_{0} (q \rarw e^{\frac{i 2 \pi}{4}} )/q^4$ associated with $K$ at $N=25000$. Real part (blue) and imaginary part (orange) of $\hat{Z}_{0}$. }
\end{center}
\end{figure}	
\noindent The q-series approaches to 
	$$ 
	\lim_{q \rarw e^{\frac{i 2 \pi}{4}}}\, \frac{2}{q^4}\hat{Z}_{0} [S^3_{-1}(K); q] \longrightarrow -4.05379317 + i 4.09952837721.
	$$
From (10), we obtain $Z_{CS} \approx 0.509582 -i 0.002858461 $.
\newline
		
\noindent At $k=5$, the analytic formula (9) gives 
	$$
	Z_{CS} \lsb S^3_{-1}(K); e^{\frac{i 2 \pi}{5}} \rsb = 0.1148764603  + i 0.3535533906. 
	$$
	We truncate the q-power series at $N=25000$ and find the limiting value of $\hat{Z}$. 
\begin{figure}[H]
\begin{center}
\includegraphics[scale=0.5]{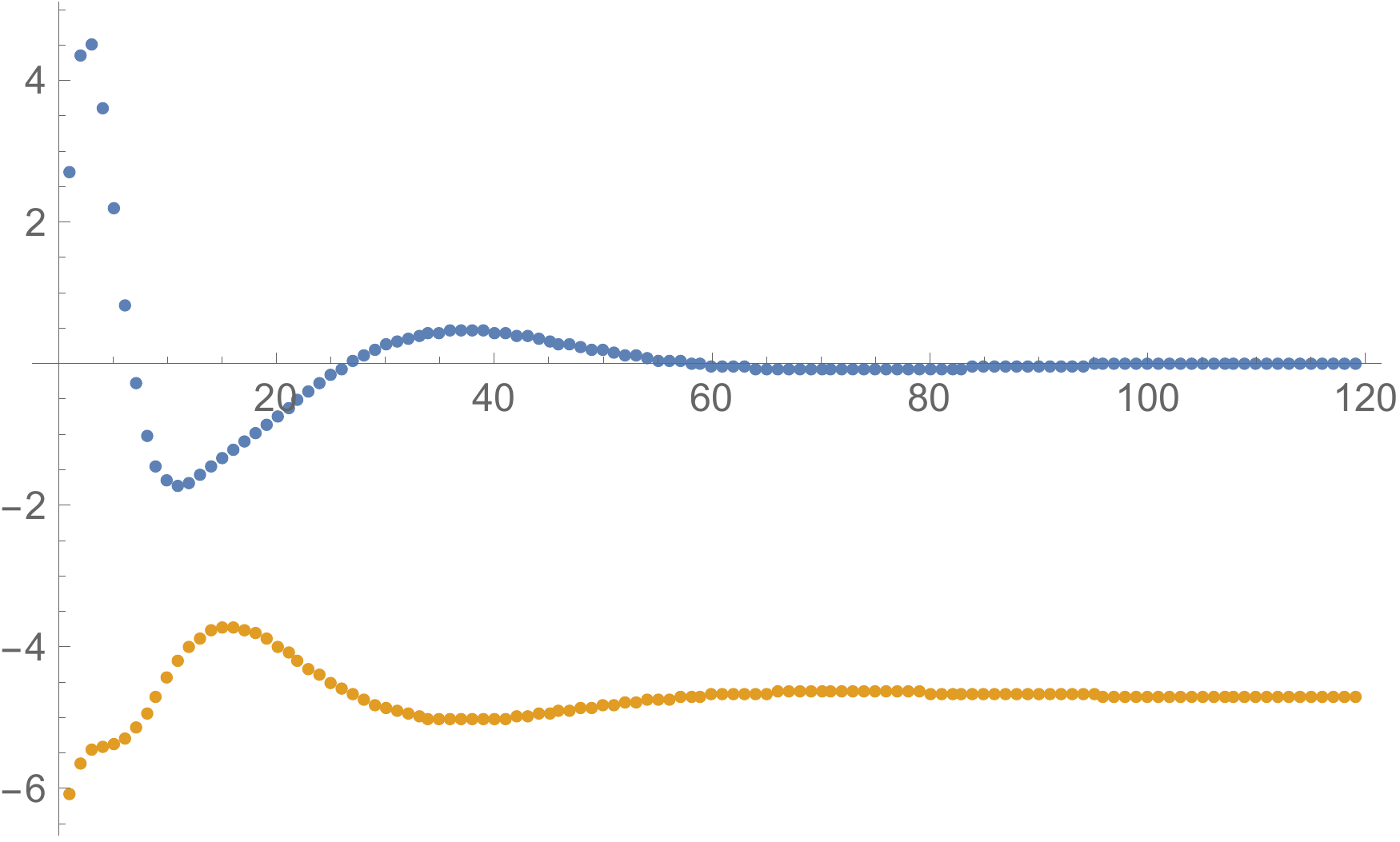}
\caption{The extrapolation of $2\hat{Z}_{0} (q \rarw e^{\frac{i 2 \pi}{5}} )/q^4 $ associated with $K$ at $N=25000$. Real part (blue) and imaginary part (orange) of $\hat{Z}_{0}$. }
\end{center}
\end{figure}	
\noindent The q-series approaches to 
	$$ 
	\lim_{q \rarw e^{\frac{i 2 \pi}{5}}}\, \frac{2}{q^4}\hat{Z}_{0} [S^3_{-1}(K); q] \longrightarrow 0.007799372126 - i 4.707580478.
	$$
From (10), $Z_{CS} \approx 0.114412 + i 0.354142$.
\newline

\appendix
\section*{Appendix}
\addcontentsline{toc}{section}{Appendix}

\section{Further plots}

We list more plots for the connected sums of knots analyzed in Section 4. In the section, the upper bound plotted on the horizontal axis correspond to the maximum value among upper bounds of summations in $F_{K_1  \, \# \, K_2} (xq^j,q),\quad j=0,\cdots ,n$ where n is an order of a $\hat{A}$-polynomial of $K_1  \, \# \, K_2$.
\begin{figure}[h!]
\begin{center}
\includegraphics[scale=1]{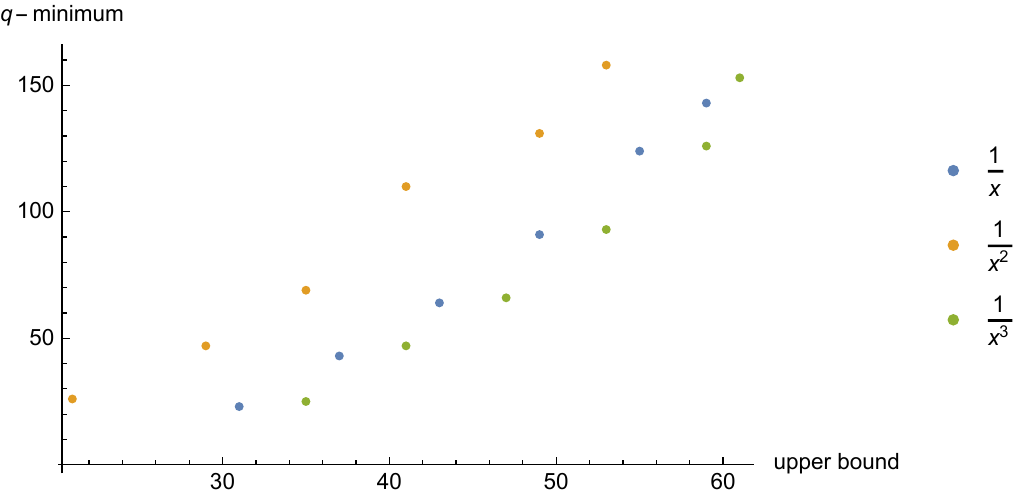}
\caption{Other powers of x in the recursion relation (5) for $K=T(2,3)\, \# \, T(2,3)$.}
\end{center}
\end{figure}

\begin{figure}[h!]
\begin{center}
\includegraphics[scale=1]{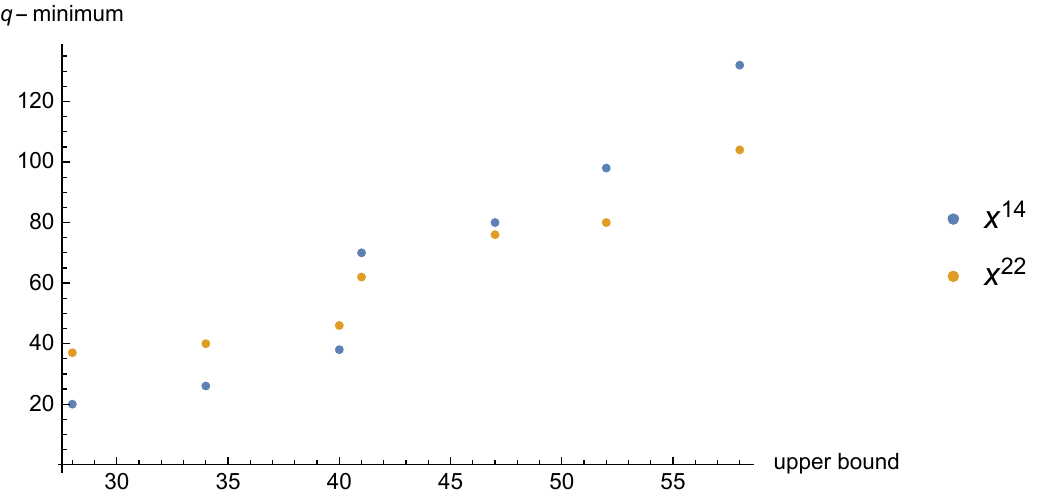}
\caption{Other powers of x in the recursion relation (5) for $K=T(2,3)\, \#\, T(2,3)$}
\end{center}
\end{figure}

\begin{figure}[h!]
\begin{center}
\includegraphics[scale=1]{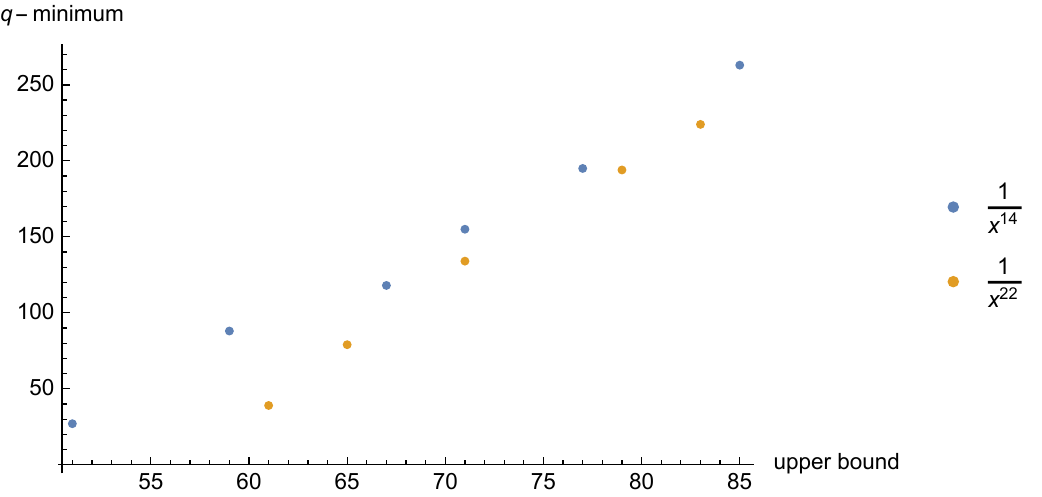}
\caption{Other powers of x in the recursion relation (5) for $K=T(2,3)\, \#\, T(2,3)$}
\end{center}
\end{figure}

\begin{figure}[h!]
\begin{center}
\includegraphics[scale=1]{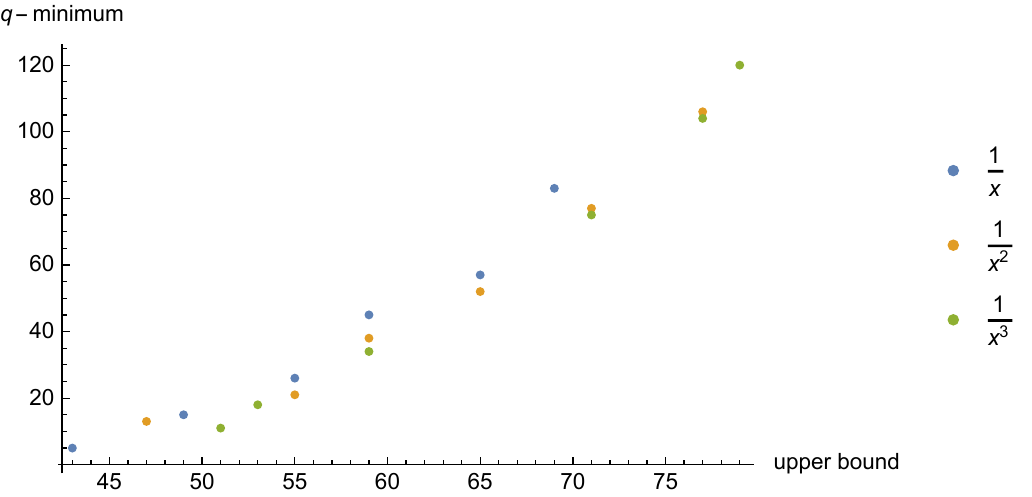}
\caption{Other powers of x in the recursion relation (6) for $K=T(2,3)\, \#\, T(2,5)$}
\end{center}
\end{figure}

\begin{figure}[h!]
\begin{center}
\includegraphics[scale=1]{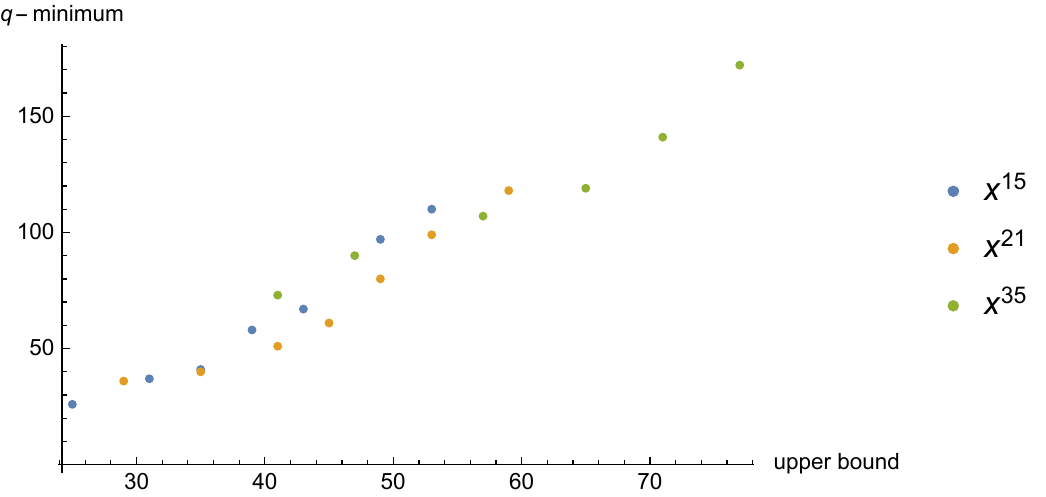}
\caption{Other powers of x in the recursion relation (6) for $K=T(2,3)\, \#\, T(2,5)$}
\end{center}
\end{figure}

\begin{figure}[h!]
\begin{center}
\includegraphics[scale=1]{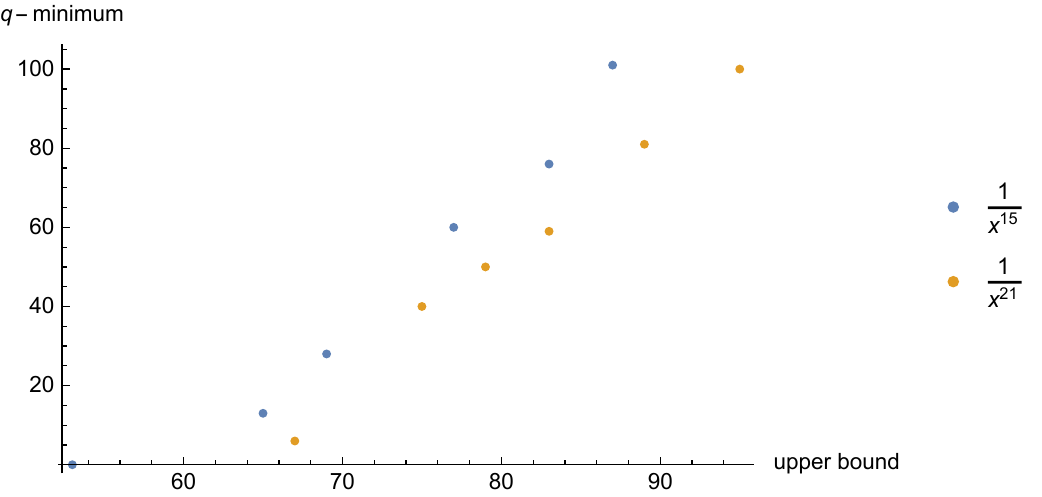}
\caption{Other powers of x in the recursion relation (6) for $K=T(2,3)\, \#\, T(2,5)$}
\end{center}
\end{figure}

\begin{figure}[h!]
\begin{center}
\includegraphics[scale=1]{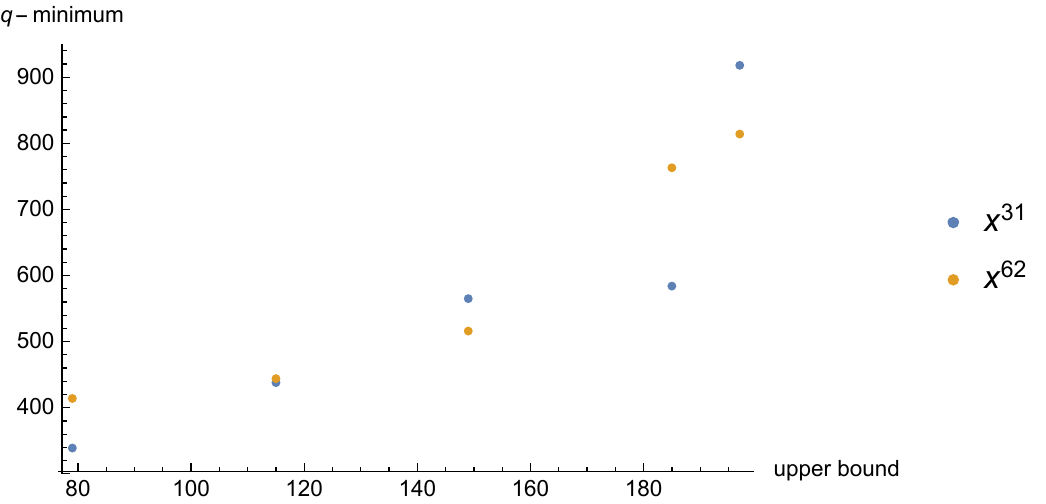}
\caption{Other powers of x in the recursion relation (7) for $K=T(2,3)\, \#\, T(3,5)$}
\end{center}
\end{figure}

\begin{figure}[h!]
\begin{center}
\includegraphics[scale=1]{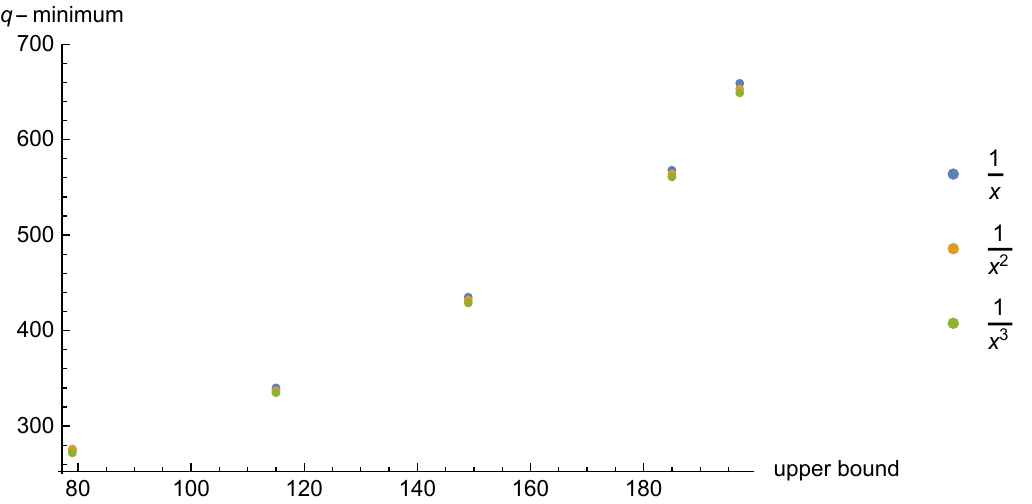}
\caption{Other powers of x in the recursion relation (7) for $K=T(2,3)\, \#\, T(3,5)$. The three dots are overlapping.}
\end{center}
\end{figure}

\begin{figure}[h!]
\begin{center}
\includegraphics[scale=1]{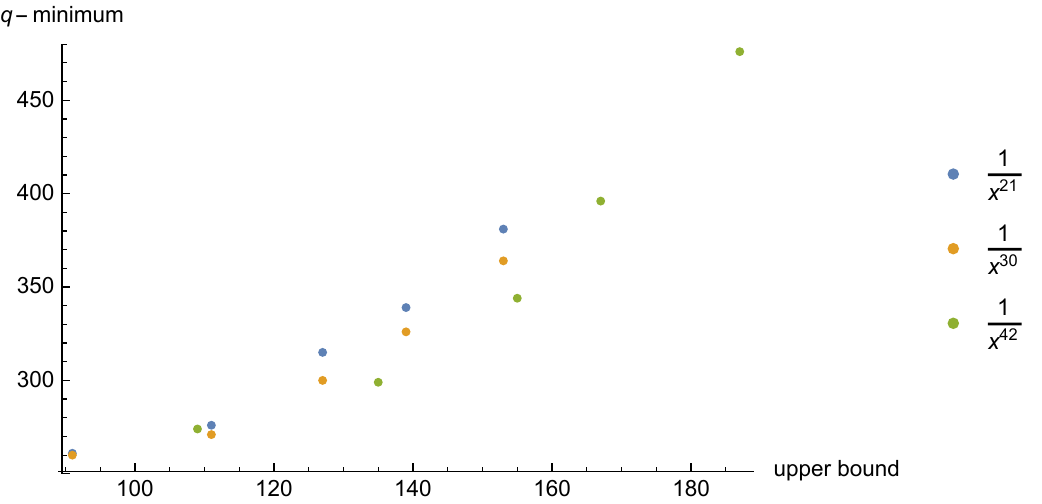}
\caption{Other powers of x in the recursion relation (7) for $K=T(2,3)\, \#\, T(3,5)$.}
\end{center}
\end{figure}

\begin{figure}[h!]
\begin{center}
\includegraphics[scale=1]{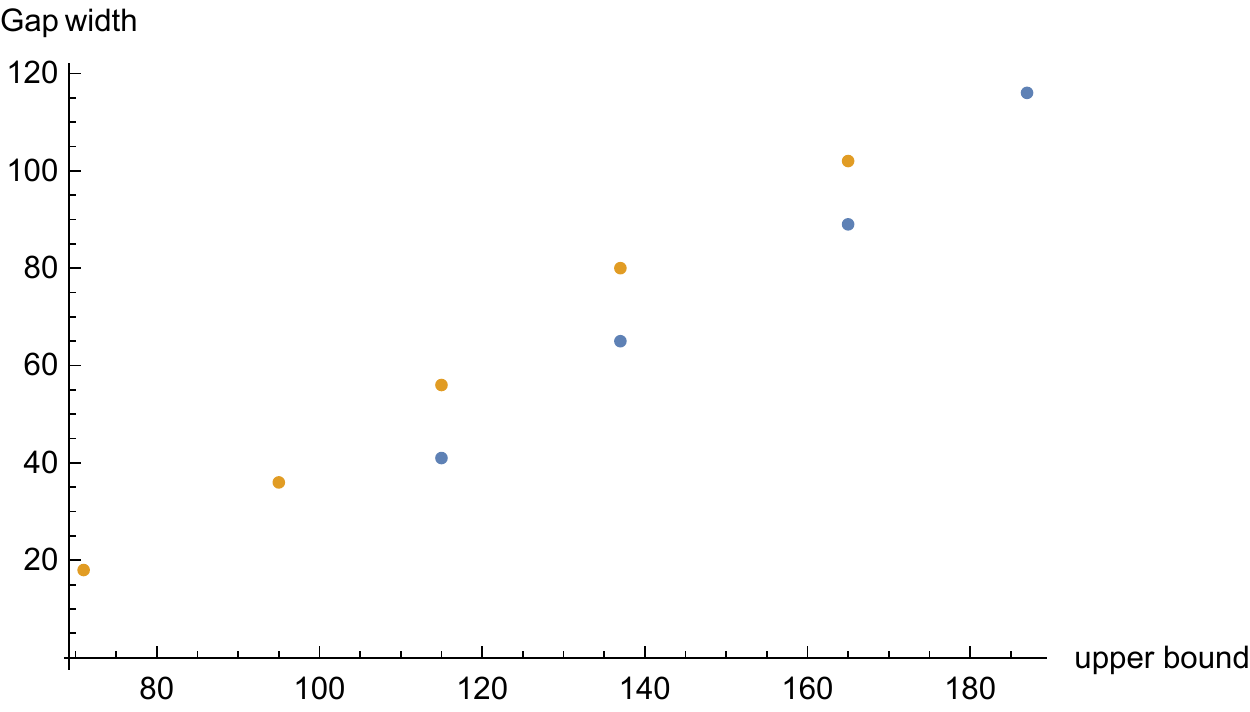}
\caption{For $K=T(2,3)\, \#\, T(2,-3)$, at $q$, the width of the gaps in $x^p$ terms (blue) and in $1/x^p$ terms (orange), $p \in \intg_{+}$ is shown. Cancellations for the blue data do not occur when the upper bounds are 71 and 95.}
\end{center}
\end{figure}

\begin{figure}[h!]
\begin{center}
\includegraphics[scale=1]{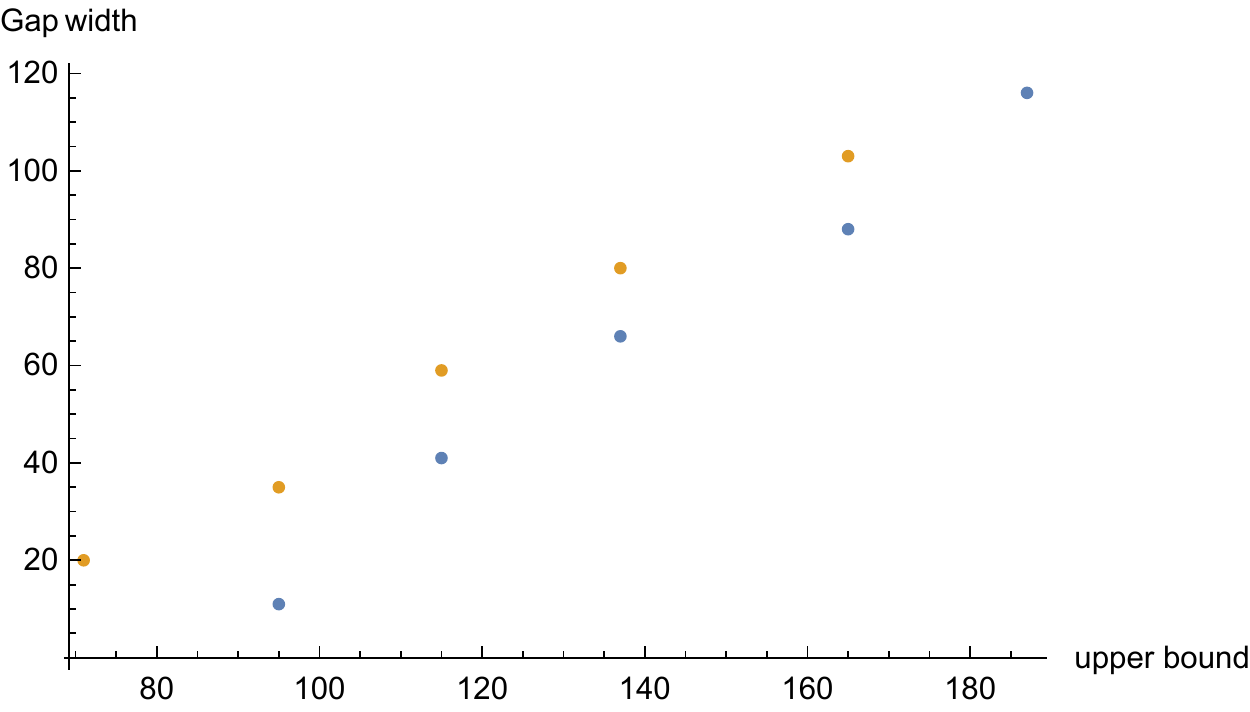}
\caption{For $K=T(2,3)\, \#\, T(2,-3)$, at $q^{2}$, the width of the gaps in $x^p$ terms (blue) and in $1/x^p$ terms (orange), $p \in \intg_{+}$ is shown. Cancellations for the blue data do not occur when the upper bound is 71.}
\end{center}
\end{figure}

\begin{figure}[h!]
\begin{center}
\includegraphics[scale=1]{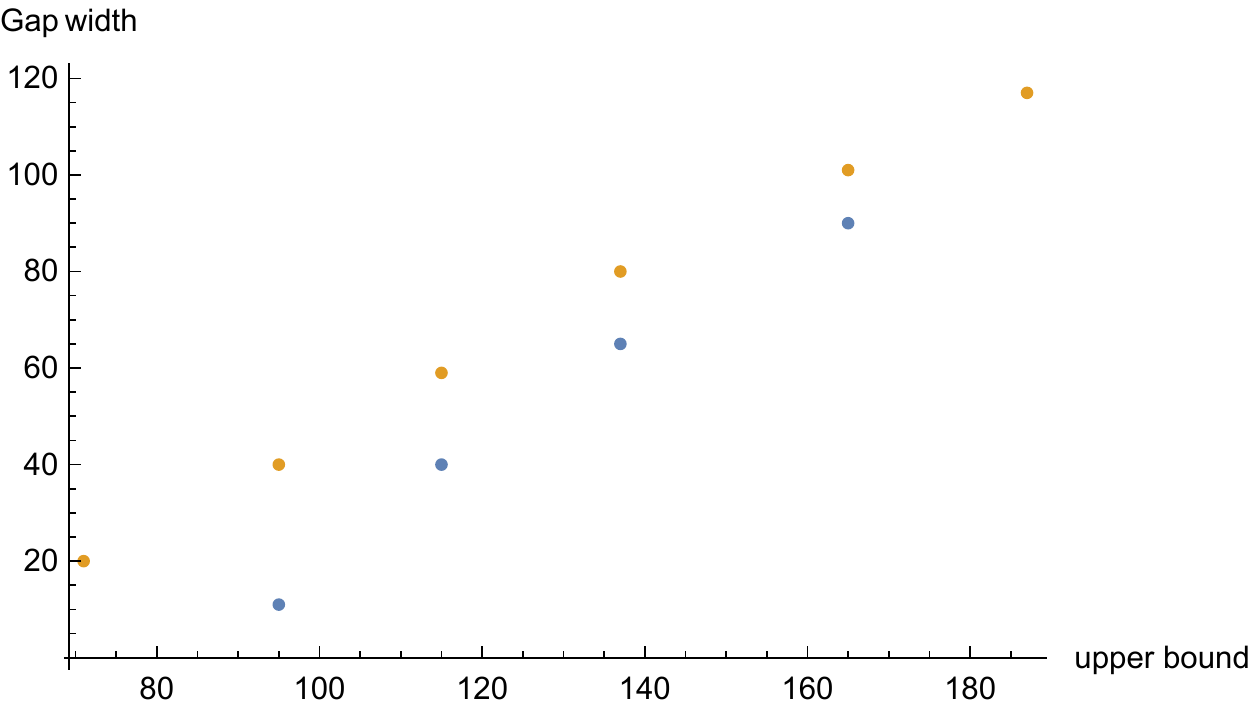}
\caption{For $K=T(2,3)\, \#\, T(2,-3)$, at $q^{3}$, the width of the gaps in $x^p$ terms (blue) and in $1/x^p$ terms (orange), $p \in \intg_{+}$ is shown. Cancellations for the blue data do not occur when the upper bound is 71.}
\end{center}
\end{figure}

\begin{figure}[h!]
\begin{center}
\includegraphics[scale=1]{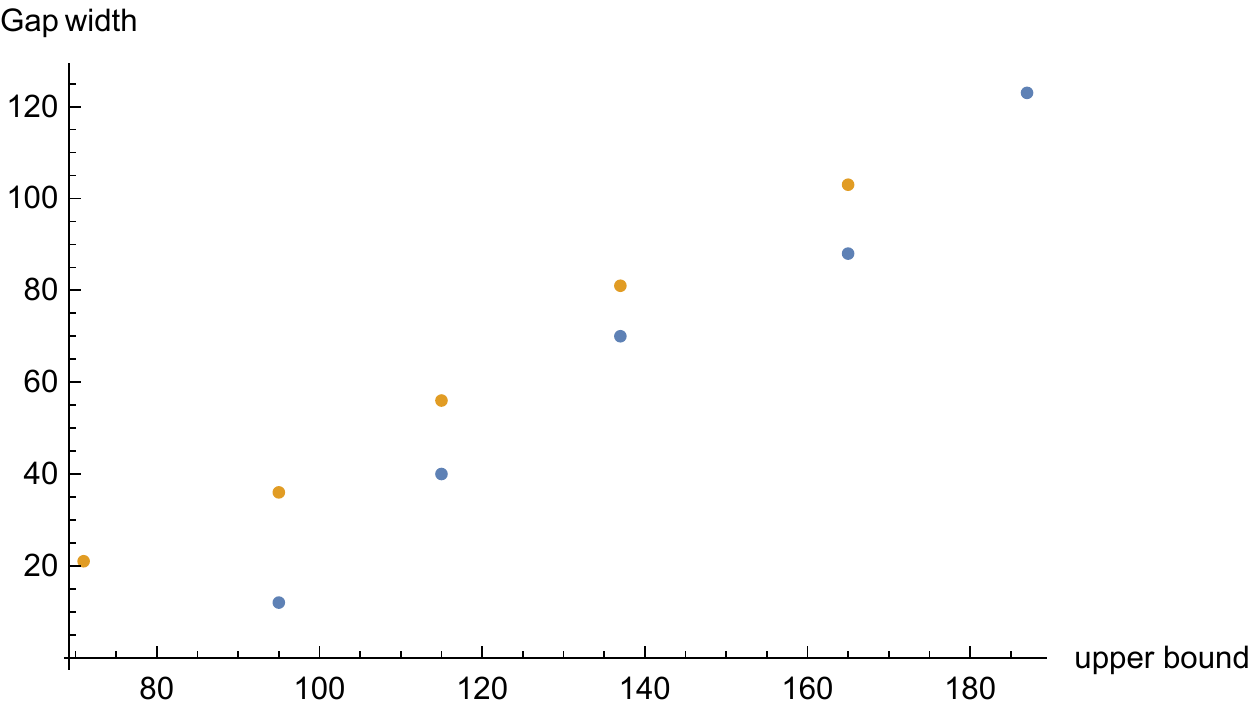}
\caption{For $K=T(2,3)\, \#\, T(2,-3)$, at $q^{-1}$, the width of the gaps in $x^p$ terms (blue) and in $1/x^p$ terms (orange), $p \in \intg_{+}$ is shown. Cancellations for the blue data do not occur when the upper bound is 71.}
\end{center}
\end{figure}

\begin{figure}[h!]
\begin{center}
\includegraphics[scale=1]{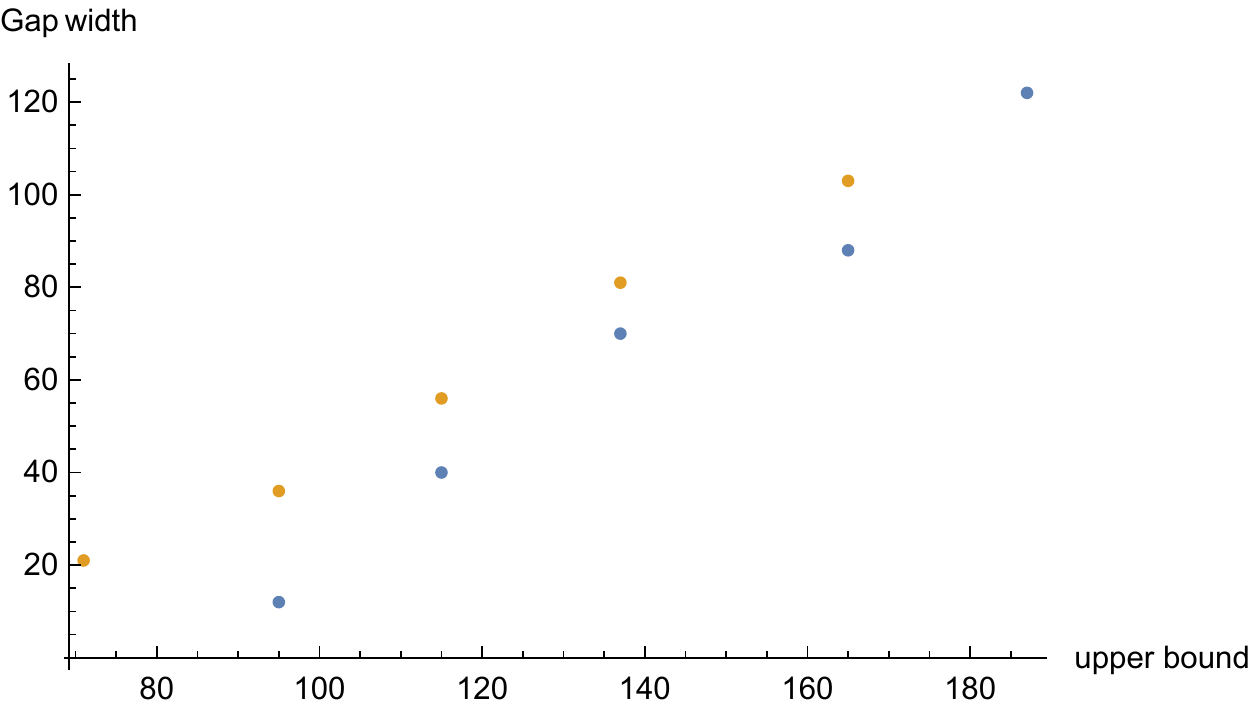}
\caption{For $K=T(2,3)\, \#\, T(2,-3)$, at $q^{-2}$, the width of the gaps in $x^p$ terms (blue) and in $1/x^p$ terms (orange), $p \in \intg_{+}$ is shown.  Cancellations for the blue data do not occur when the upper bound is 71.}
\end{center}
\end{figure}

\begin{figure}[h!]
\begin{center}
\includegraphics[scale=1]{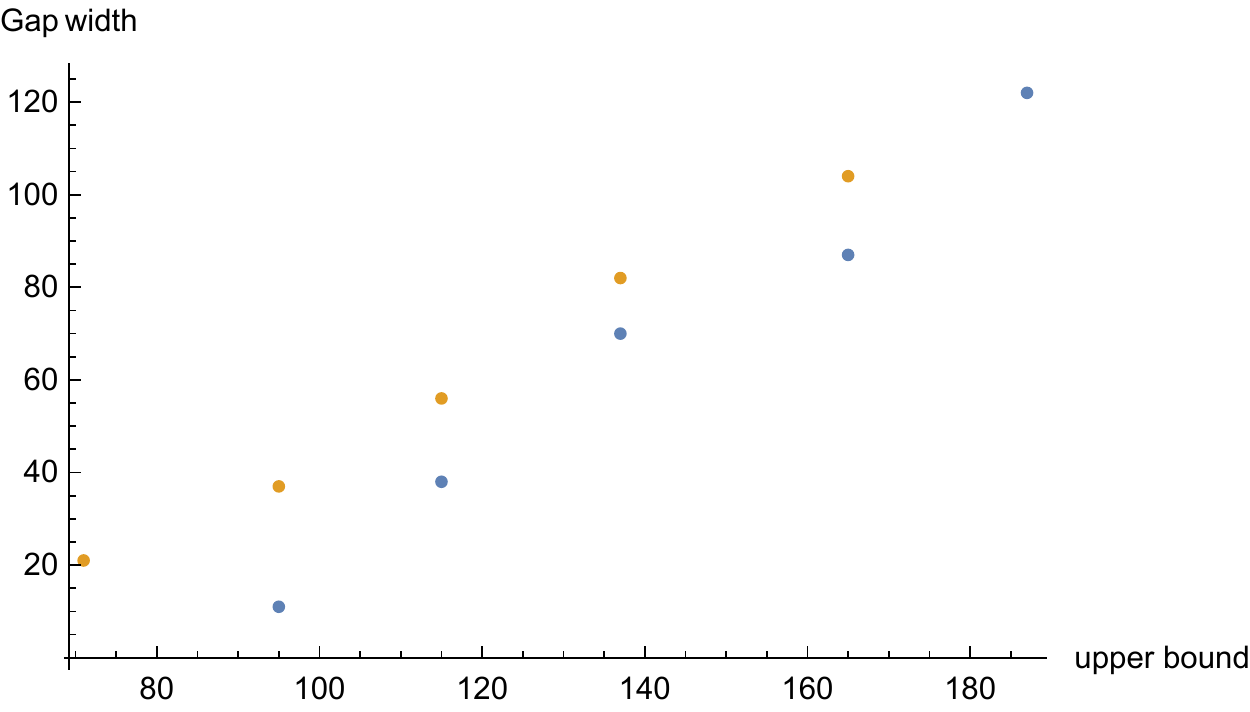}
\caption{For $K=T(2,3)\, \#\, T(2,-3)$, at $q^{-3}$, the width of the gaps in $x^p$ terms (blue) and in $1/x^p$ terms (orange), $p \in \intg_{+}$ is shown. Cancellations for the blue data do not occur when the upper bound is 71.}
\end{center}
\end{figure}

\begin{figure}[h!]
\begin{center}
\includegraphics[scale=1]{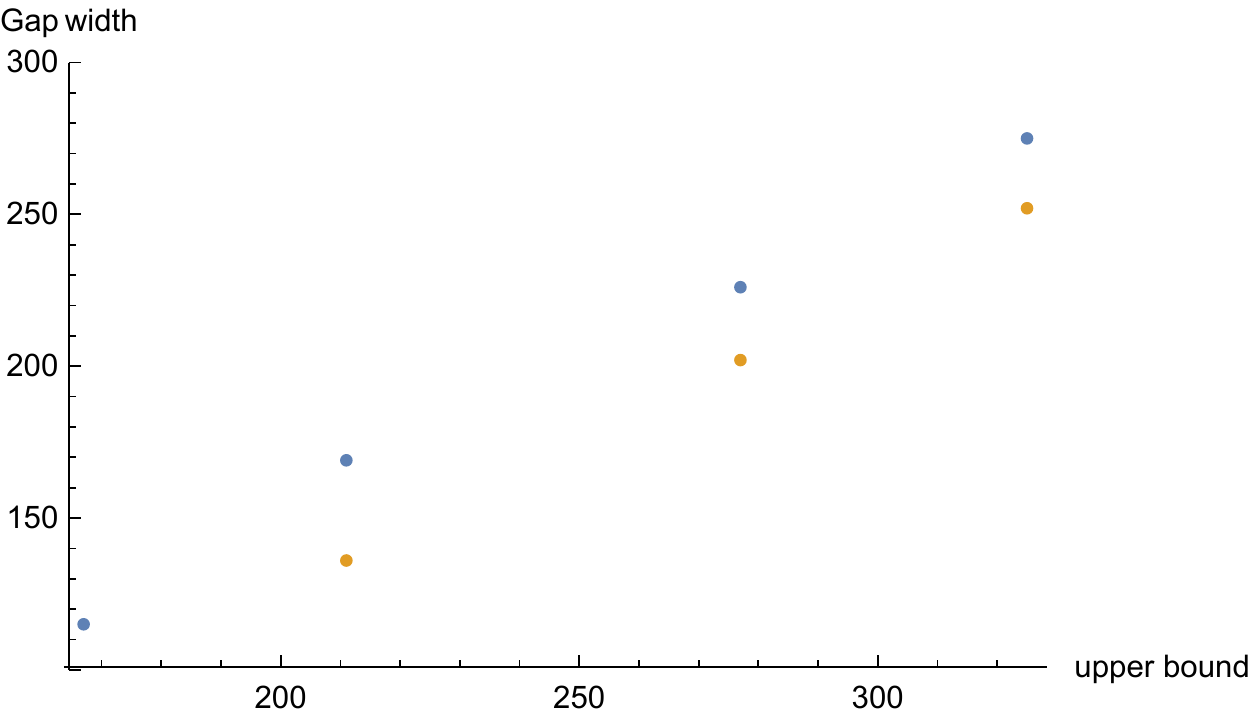}
\caption{For $K=T(2,3)\, \#\, T(2,-3)$, at $q^{300}$, the width of the gaps in $x^p$ terms (blue) and in $1/x^p$ terms (orange), $p \in \intg_{+}$ is shown.}
\end{center}
\end{figure}

\begin{figure}[h!]
\begin{center}
\includegraphics[scale=1]{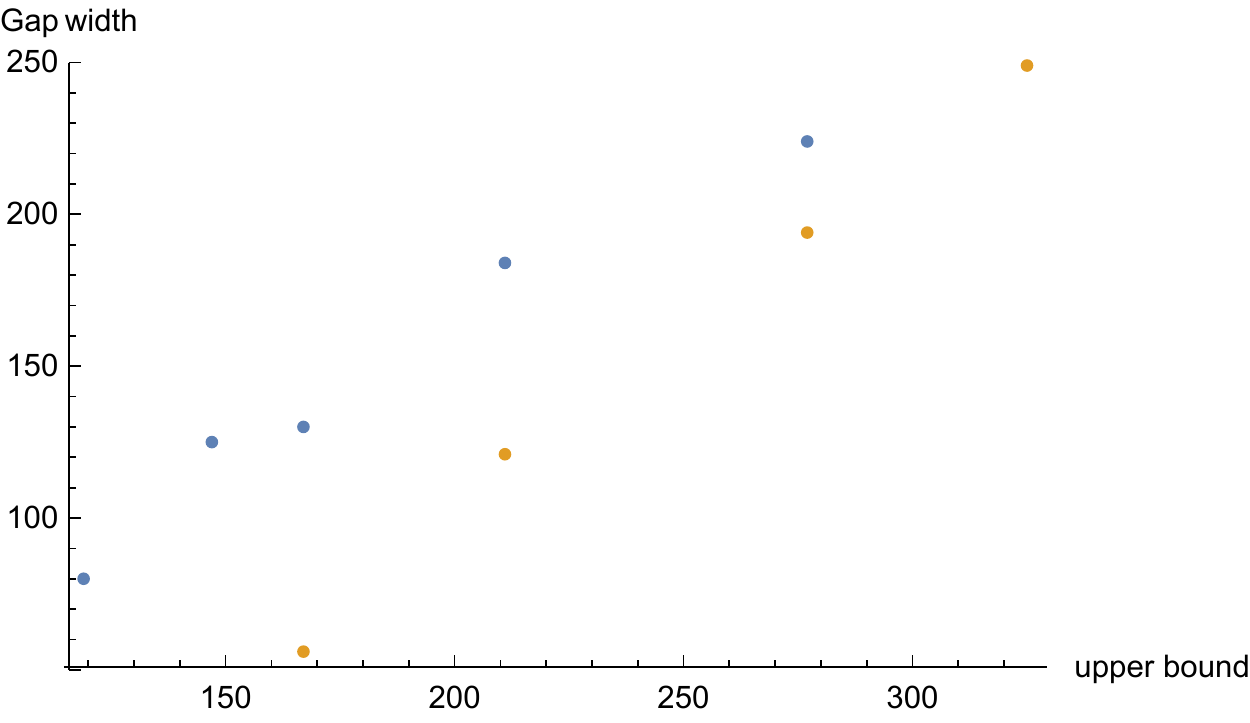}
\caption{For $K=T(2,3)\, \#\, T(2,-3)$, at $q^{500}$, the width of the gaps in $x^p$ terms (blue) and in $1/x^p$ terms (orange), $p \in \intg_{+}$ is displayed. Cancellations for the orange data do not occur when upper bounds are 95 and 121.}
\end{center}
\end{figure}

\begin{figure}[h!]
\begin{center}
\includegraphics[scale=1]{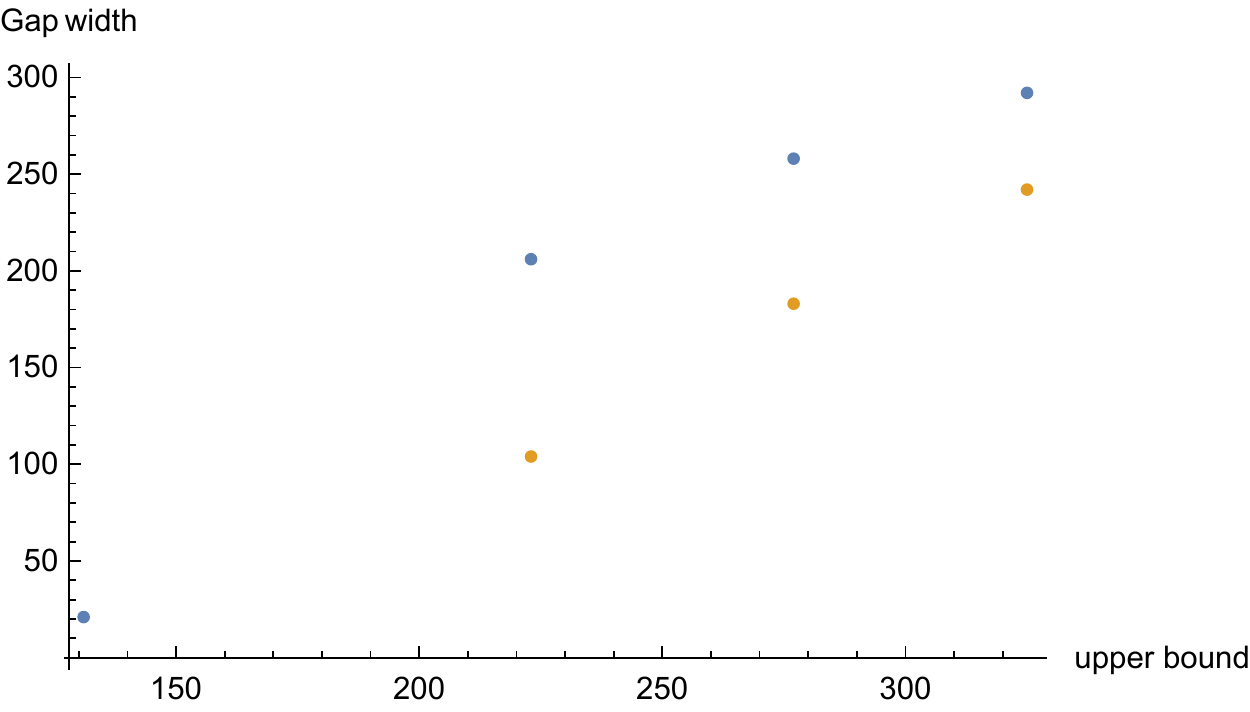}
\caption{For $K=T(2,3)\, \#\, T(2,-3)$, at $q^{793}$, the width of the gaps in $x^p$ terms (blue) and in $1/x^p$ terms (orange), $p \in \intg_{+}$ is displayed. Cancellations for the orange data do not occur when the upper bound is 107.}
\end{center}
\end{figure}

\begin{figure}[h!]
\begin{center}
\includegraphics[scale=1]{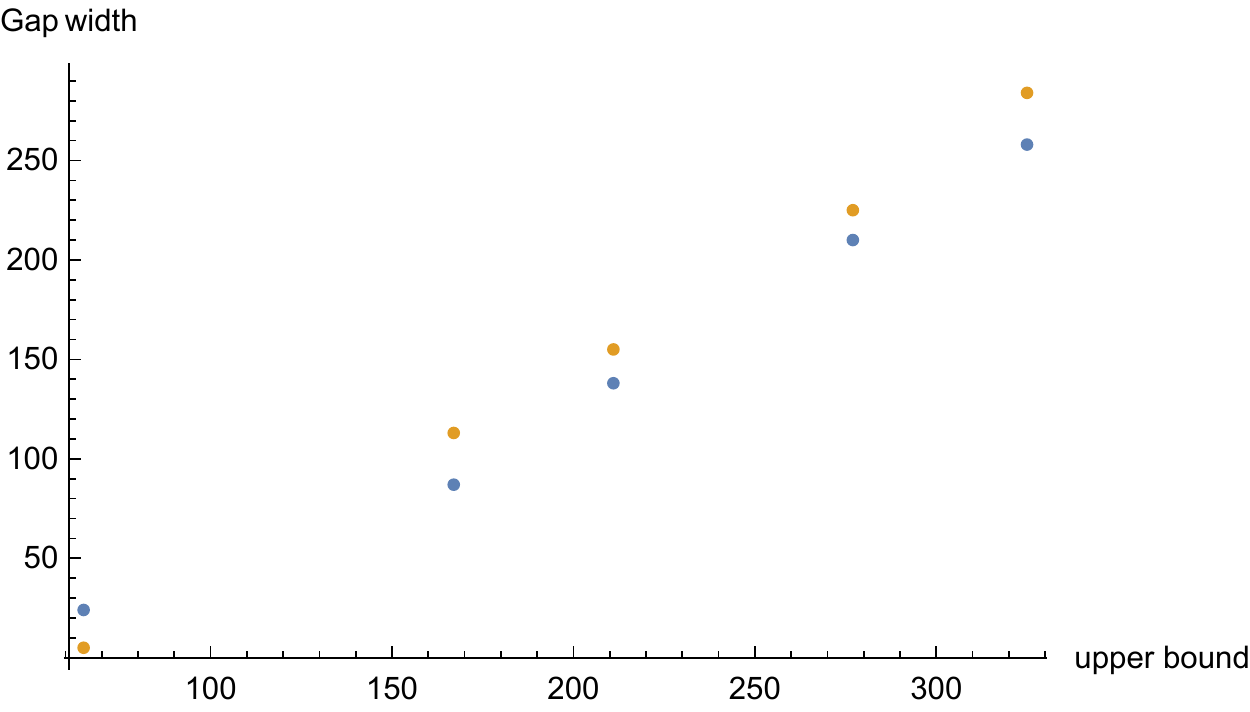}
\caption{For $K=T(2,3)\, \#\, T(2,-3)$, at $q^{-105}$, the width of the gaps in $x^p$ terms (blue) and in $1/x^p$ terms (orange), $p \in \intg_{+}$ is displayed.}
\end{center}
\end{figure}

\begin{figure}[h!]
\begin{center}
\includegraphics[scale=1]{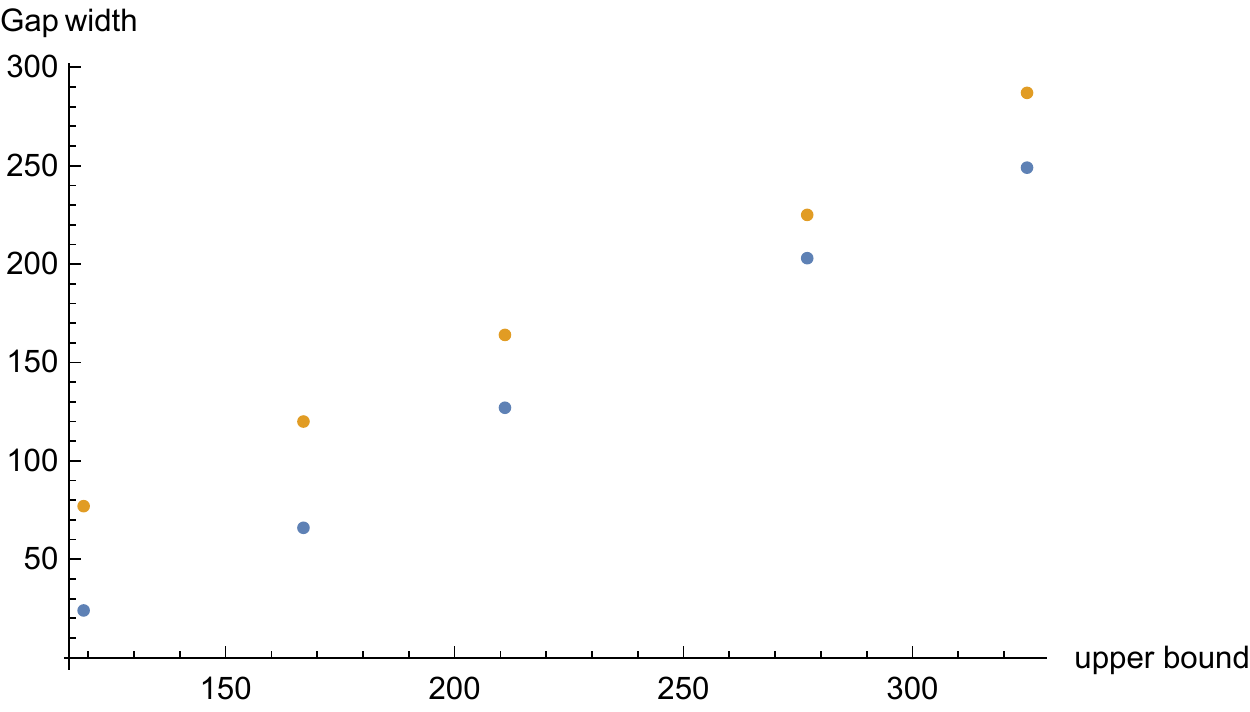}
\caption{For $K=T(2,3)\, \#\, T(2,-3)$, at $q^{-300}$, the width of the gaps in $x^p$ terms (blue) and in $1/x^p$ terms (orange), $p \in \intg_{+}$ is displayed.}
\end{center}
\end{figure}

\begin{figure}[h!]
\begin{center}
\includegraphics[scale=1]{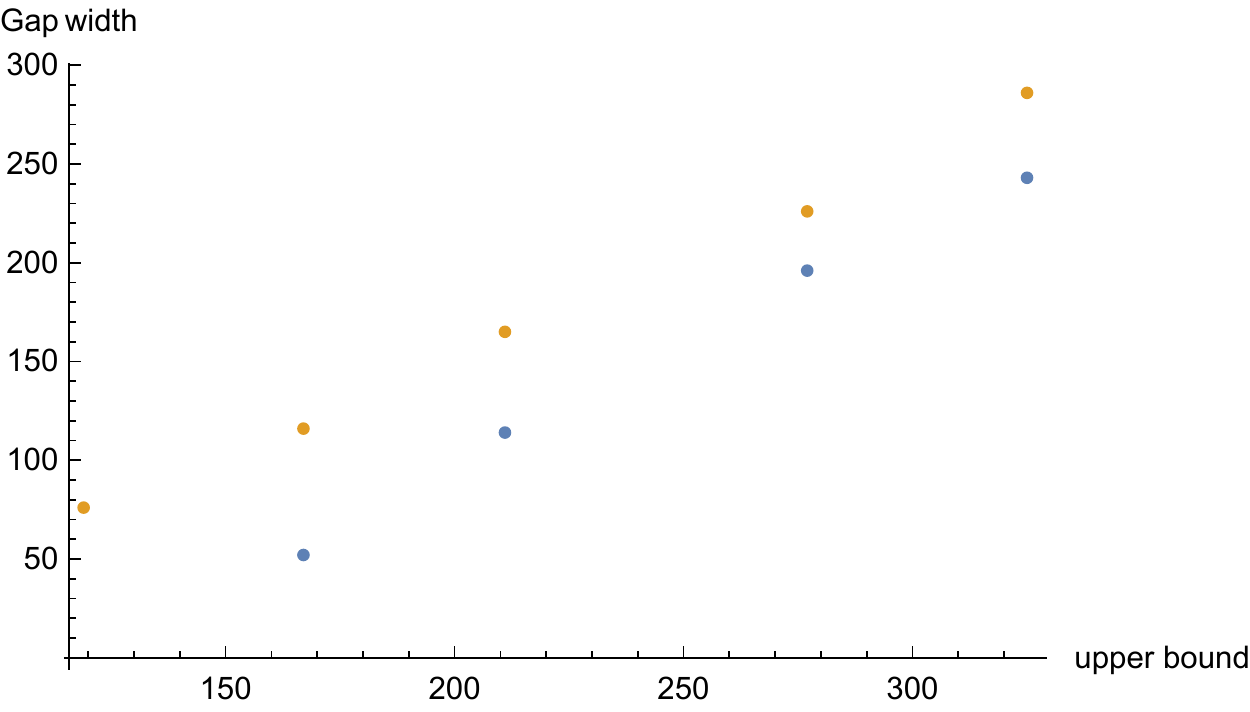}
\caption{For $K=T(2,3)\, \#\, T(2,-3)$, at $q^{-400}$, the width of the gaps in $x^p$ terms (blue) and in $1/x^p$ terms (orange), $p \in \intg_{+}$ is displayed. Cancellations for the blue data do not occur when the upper bound is 95.}
\end{center}
\end{figure}
\clearpage

The five terms in (8) at $q^{100}$ are recorded below. Due to their lengthy expressions, ellipsis are used.
\begin{align*}
H_0 & = -\frac{5}{x^{269}}+\frac{33}{x^{268}}-\frac{61}{x^{267}}+ \cdots +\frac{167}{x^3}-\frac{89}{x^2}-\frac{99}{x}+161 -62x-63 x^2+140 x^3-111 x^4+43 x^5\\
& -5 x^6-64 x^7+195 x^8-335 x^9+375 x^{10}-314 x^{11}+146 x^{12}+114x^{13}-295 x^{14}+415 x^{15}\\
& -461 x^{16}+ \cdots +33 x^{284}-61 x^{285}+38 x^{286}
\end{align*}

\begin{align*}
H_1 & = \frac{77}{x^{269}}-\frac{81}{x^{267}}+\frac{9}{x^{266}}+ \cdots +\frac{238}{x^3}-\frac{52}{x^2}-\frac{248}{x}+429 -276 x-200 x^2+444 x^3-212 x^4-124 x^5\\
& +451 x^6 -101 x^7-812 x^8+940 x^9-489 x^{10}+369 x^{11}-239x^{12}-224 x^{13}+598 x^{14}-626 x^{15}+350 x^{16}\\
& + \cdots +7 x^{284}-84 x^{285}+7 x^{286}
\end{align*}

\begin{align*}
H_2 & = -\frac{33}{x^{269}}-\frac{131}{x^{268}}+\frac{61}{x^{267}} + \cdots -\frac{228}{x^3}-\frac{373}{x^2}+\frac{349}{x}-47 -181 x-57 x^2+275 x^3+374 x^4-386 x^5\\
& -435 x^6+367 x^7+622 x^8-1080 x^9+446 x^{10}-133 x^{11}+193x^{12}+293 x^{13}-552 x^{14}+57 x^{15}\\
& -183 x^{16}+ \cdots -128 x^{284}+63 x^{285}+92 x^{286}
\end{align*}

\begin{align*}
H_3 & = -\frac{59}{x^{269}}+\frac{65}{x^{268}}+\frac{114}{x^{267}} + \cdots -\frac{167}{x^3}+\frac{374}{x^2}+\frac{50}{x}-429 +367 x+229 x^2-682 x^3+62 x^4+363 x^5\\
& +81 x^6-489 x^7-105 x^8+1016 x^9-798 x^{10}+383 x^{11}-216x^{12}-313 x^{13}+557 x^{14}-319 x^{15}\\
& +707 x^{16} + \cdots -84 x^{283}+56 x^{284}+115 x^{285}-122 x^{286}
\end{align*}

\begin{align*}
H_4 & = \frac{20}{x^{269}}+\frac{33}{x^{268}}-\frac{33}{x^{267}} + \cdots -\frac{10}{x^3}+\frac{140}{x^2}-\frac{52}{x}-114 +152 x+89 x^2-185 x^3-97 x^4+195 x^5\\
& -175 x^6+136 x^7 +169 x^8-387 x^9+395 x^{10}-320 x^{11}+115 x^{12}+131 x^{13}-308 x^{14}+473 x^{15}\\
& -413 x^{16} + \cdots +17 x^{283}+32 x^{284}-33 x^{285}-15 x^{286}
\end{align*}

\end{document}